\documentclass[pdflatex,sn-mathphys-num]{sn-jnl}

\usepackage{graphicx}%
\usepackage{multirow}%
\usepackage{amsmath,amssymb,amsfonts}%
\usepackage{amsthm}%
\usepackage{mathrsfs}%
\usepackage[title]{appendix}%
\usepackage{xcolor}%
\usepackage{textcomp}%
\usepackage{manyfoot}%
\usepackage{booktabs}%
\usepackage{algorithm}%
\usepackage{algorithmicx}%
\usepackage{algpseudocode}%
\usepackage{listings}%

\usepackage{bm}

\usepackage{float} 

\begin{document}

\title[Modeling Nonlinear Dynamics from Videos]{Modeling Nonlinear Dynamics from Videos}


\author[1,2]{\fnm{Antony} \sur{Yang}}\email{thy25@cam.ac.uk}
\author[1]{\fnm{Joar} \sur{Ax\aa s}}\email{jgoeransson@ethz.ch}
\author[3,4]{\fnm{Fanni} \sur{K\'{a}d\'{a}r}}\email{fanni.kadar@mm.bme.hu}
\author[3,4]{\fnm{G\'{a}bor} \sur{St\'{e}p\'{a}n}}\email{stepan@mm.bme.hu}
\author*[1]{\fnm{George} \sur{Haller}}\email{georgehaller@ethz.ch}

\affil[1]{\orgdiv{Institute for Mechanical Systems}, \orgname{ETH Z\"{u}rich}, \orgaddress{\street{Leonhardstrasse 21}, \city{Z\"{u}rich}, \postcode{8092}, \country{Switzerland}}}

\affil[2]{\orgdiv{Cavendish Laboratory, Department of Physics}, \orgname{University of Cambridge}, \orgaddress{\city{Cambridge}, \postcode{CB2 1TN}, \country{United Kingdom}}}

\affil[3]{\orgdiv{Department of Applied Mechanics, Faculty of Mechanical Engineering}, \orgname{Budapest University of Technology and Economics}, \orgaddress{\street{Műegyetem rkp. 3}, \city{Budapest}, \postcode{H-1111}, \country{Hungary}}}

\affil[4]{\orgdiv{MTA-BME Lendület Momentum Global Dynamics Research Group}, \orgname{Budapest University of Technology and Economics}, \orgaddress{\street{Műegyetem rkp. 3}, \city{Budapest}, \postcode{H-1111}, \country{Hungary}}}


\abstract{We introduce a method for constructing reduced-order models directly from videos of dynamical systems. The method uses a non-intrusive tracking to isolate the motion of a user-selected part in the video of an autonomous dynamical system. In the space of delayed observations of this motion, we reconstruct a low-dimensional attracting spectral submanifold (SSM) whose internal dynamics serves as a mathematically justified reduced-order model for nearby motions of the full system. We obtain this model in a simple polynomial form that allows explicit identification of important physical system parameters, such as natural frequencies, linear and nonlinear damping and nonlinear stiffness. Beyond faithfully reproducing attracting steady states and limit cycles, our SSM-reduced models can also uncover  hidden motion not seen in the video, such as unstable fixed points and unstable limit cycles forming basin boundaries. We demonstrate all these features on experimental videos of five physical systems: a double pendulum, an inverted flag in counter-flow, water sloshing in tank, a wing exhibiting aeroelastic flutter and  a shimmying wheel.}




\maketitle

\section{Introduction}\label{sec1}

Nonlinear dynamical systems are prevalent in numerous fields of nature and engineering.
Examples include sloshing in tank trucks \cite{truck_slosh}, turbulent flows \cite{intro_turbulence}, spacecraft motion \cite{intro_spacecraft}, and vibrations in jointed structures \cite{intro_bolt_joints}.
Full-order modeling of such systems is often challenging due to their high number of degrees of freedom and uncertain physical parameters. For these reasons, reduced-order modeling of nonlinear mechanical systems has been an increasingly active area of research that promises major benefits in system identification \cite{breunung23}, design optimization \cite{detroux21} and model-predictive control \cite{alora23b}. A further recent boost to this effort is the trend is to construct digital twins, i.e., highly accurate, interpretable and predictive models of the current states of physical assets \cite{perno22}.

A common approach to reducing nonlinear dynamical systems to lower-dimensional models is the proper orthogonal decomposition (POD) followed by a Galerkin projection \cite{intro_pod1,intro_pod2}. This approach implicitly assumes   that model subspaces of the linearized system remain nearly invariant under inclusion of nonlinearities, which is a priori unknown and can often be secured only by selecting an unnecessarily large number of linear modes. Another popular model reduction technique is the Dynamic Mode Decomposition (DMD) and its variants \cite{DMD,schmid22}, which find the best fitting linear autonomous dynamical system for the available observable data. Passing to the dominant modes of this linear system then provides a linear reduced-order model of the full dynamics. While efficient in capturing linearizable dynamics, DMD methods return linear systems and hence cannot model intrinsically nonlinear phenomena, such as coexisting isolated attracting fixed points, limit cycles or transitions between such states \cite{koopman_fail,liu23}.

Machine learning methods based on the training of neural networks \cite{machine_learn_model,machine_learn_model2} have also been explored as data-driven model reduction alternatives, but they often lack physical interpretability, are prone to overfitting and require large amounts of data and extensive tuning. Within the category of machine learning, sparse identification of nonlinear dynamics (SINDy) \cite{SINDy} can provide interpretable models but only if an appropriate reduced set of variables is already known. Even in that case, however, the outcome of the process is generally sensitive to the choice of sparsification parameters.

In recent years, reduction to spectral submanifolds (SSMs) has appeared as an alternative for constructing reduced-order models for nonlinear dynamical systems from data \cite{ssm_data_old,SSMLearn,SSMLearn_reduction,fastSSM}.
A primary SSM of a dynamical system is the unique smoothest invariant manifold tangent to a nonresonant spectral subspace of the linearized system at a steady state. Such manifolds have long been envisioned and formally approximated as nonlinear normal modes (NNMs) in a series of papers initiated by the seminal work of Shaw and Pierre \cite{shaw93} in the 1990's. The existence, uniqueness and smoothness of such manifolds, however,  has only been clarified more recently for various types of steady states \cite{cabre03,haro06,ssm,ssm_mixed} and general external forcing \cite{haller23b}. Importantly, the internal dynamics of attracting SSMs tangent to a span of slowest eigenmodes provide a mathematically exact reduced-order model with which all nearby trajectories synchronize exponentially fast.

Analytic and data-driven tools, available as open-source MATLAB codes, have been developed for constructing SSM-based reduced-order models. For analytic treatments, the \emph{SSMTool} \cite{SSMTool2} package computes SSMs directly from governing equations, whereas the data-driven methods \emph{SSMLearn} \cite{SSMLearn} and \emph{fastSSM} \cite{fastSSM} construct models from time-series data. These SSM-based reduced order modeling algorithms have been successfully applied to several nonlinear systems known from numerical or experimental data.
Examples include the flow past a cylinder \cite{SSMLearn}, a Brake-Reuss beam \cite{SSMLearn_reduction}, water tank sloshing \cite{fastSSM}, transition to turbulence in pipe flows \cite{kaszas24} and finite element models of various structures going beyond a million degrees of freedom  \cite{SSMTool2,cenedese24}.

While all these studies demonstrate the applicability and robustness of SSM models inferred from data, the scarce availability of experimental measurements many physical settings motivates the extension of SSM-based model reduction to general video data. Such an extension should ideally work with general video footage, such as an excerpt from a documentary or instructional video, that was not necessarily generated in a sterile environment for the sole purpose of model reduction. An additional benefit of a purely video-based model reduction would be its nonintrusive nature. This is especially important for slender structures where an attached sensor would alter the behavior of the system. In other cases, video-based system modeling may be the only viable option because  the placement of a reliable physical sensor is unrealistic due to extreme temperatures, pressure or humidity.

While such techniques have been widely used in various engineering and scientific disciplines, such as robotics \cite{cv_robotics}, medical imaging \cite{cv_medical}, manufacturing \cite{cv_manufacturing}, autonomous driving \cite{cv_autodrive}, structural health monitoring \cite{SHM1}, fluid mechanics (particle image velocimetry) \cite{slosh_experiment}, and unsupervised physical scene understanding \cite{physics_from_video1}, their application in reduced-order modeling of nonlinear dynamical systems has remained largely unexplored. Indeed, video-based modeling has been mostly tried for identifying parameters in systems whose governing equations are already known (see \cite{nonlinear_dynamics_from_video} for a review). The same reference  proposes a method to infer nonlinear dynamics by training neural networks on  videos generated by solutions of simple, polynomial ODEs \cite{nonlinear_dynamics_from_video}. So far, however, no applications of these methods to videos of  real physical experiments have appeared. 

A major challenge in video-based reduced-order modeling modeling is the visual tracking of physical systems, which is complicated by changes in object motion, cluttered backgrounds, occlusion, and changes in target appearance \cite{review_correlation_track,review_vision_track,review_deep_learning_track}.
Generic tracking, also known as short-term or model-free tracking, involves continuously localizing a target in a video sequence using a single example of its appearance, usually initialized in the first frame.
While nearly all existing tracking algorithms focus on object localization, our goal here is to use visual tracking as an experimental sensing method for dynamical systems.
In other words, we are interested in the accurate recovery of trajectories rather than locating objects in videos.

\begin{figure}[t]
    \centering
    \includegraphics[width=1\linewidth]{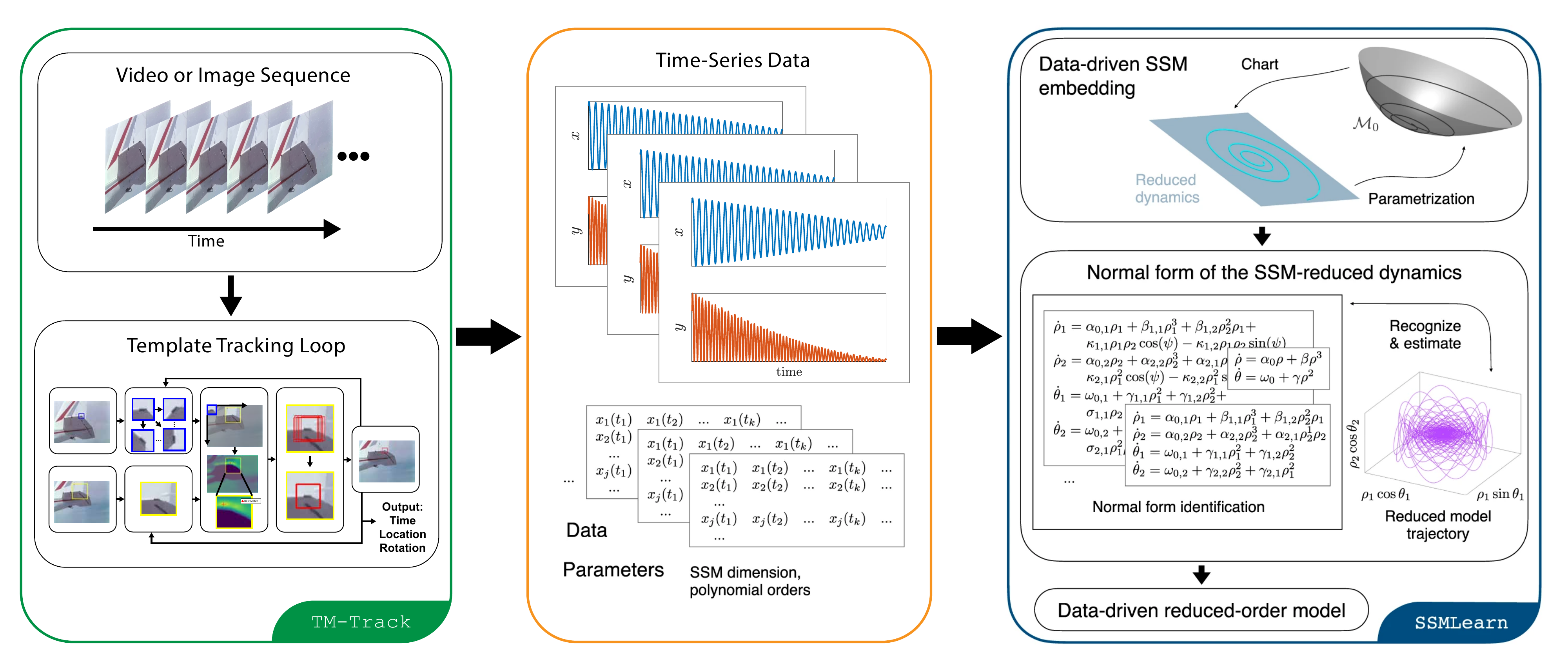}
    \caption[SSMLearn]%
    {\textbf{Schematic of our Tracking and Modeling Approach} --- (green) First, we prepare input videos or image sequences that capture transient phenomena. We select a template on the video to track.
    The tracking algorithm outputs the $(x,y)$ pixel and template rotation as time-series data
    (orange) Next,  the data is pre-processed with possible noise reduction and delay-embedded to produce a suitable observable space.
    (blue) Finally, we train a reduced-order model with the \emph{SSMLearn} algorithm \cite{SSMLearn} in the delay-embedded space to make predictions for previously unseen initial conditions or to predict behavior under additional external forcing.}
    \label{fig:ssmlearn}
\end{figure}

The tracking algorithms developed in recent decades, such as correlation-based methods \cite{kcf,circular_kcf,MOSSE,dcf_scale,CSRT} and convolutional neural network-based (CNN) methods \cite{track_CNN_offline,track_CNN_offline2,track_CNN_online,track_CNN_multi,RCNN}, are founded on the assumption that the appearance of the tracked object changes over time, requiring a process that learns such changes.
The main objective of these methods is to continuously update the model of the tracked object, learning its evolving shape and size, to ensure its robust and accurate localization. 
However, due to the changing nature of the object's representation, the tracked features are not constant observables on the phase space, and are therefore typically unsuitable for system identification. In response to the challenge of precise trajectory extraction from videos for dynamical problems, Kara et al. \cite{granular_track} proposed a tracking algorithm based on deep learning, addressing the identity switching problem encountered in particle tracking scenarios such as walking droplets and granular intruders experiments. Their approach, however, relies on a machine learning model that necessitates manual labeling of data for training the tracker.

In this work, we introduce a tracking algorithm that prioritizes trajectory recovery to assemble data for SSM-based model reduction from videos of physical systems. The tracking method is based on the classic template matching technique \cite{computer_vision}, which achieves marker-free and rotation-aware tracking that requires only a one-time initialization from the user. With the trajectories obtained from our proposed tracker, we train SSM-reduced polynomial models, employing the open-source MATLAB package \emph{SSMLearn}. The main steps of this procedure are summarized in Figure \ref{fig:ssmlearn}.
We illustrate on several examples that the methodology developed here is applicable to videos of a diverse set of physical systems that ranging from a double pendulum and an inverted flag through water sloshing and aeroelastic flutter to wheel shimmy.

\section{Marker-free tracking with template matching}\label{section:markerless_track}
\begin{figure}[h]
    \centering
    \includegraphics[width=1\linewidth]{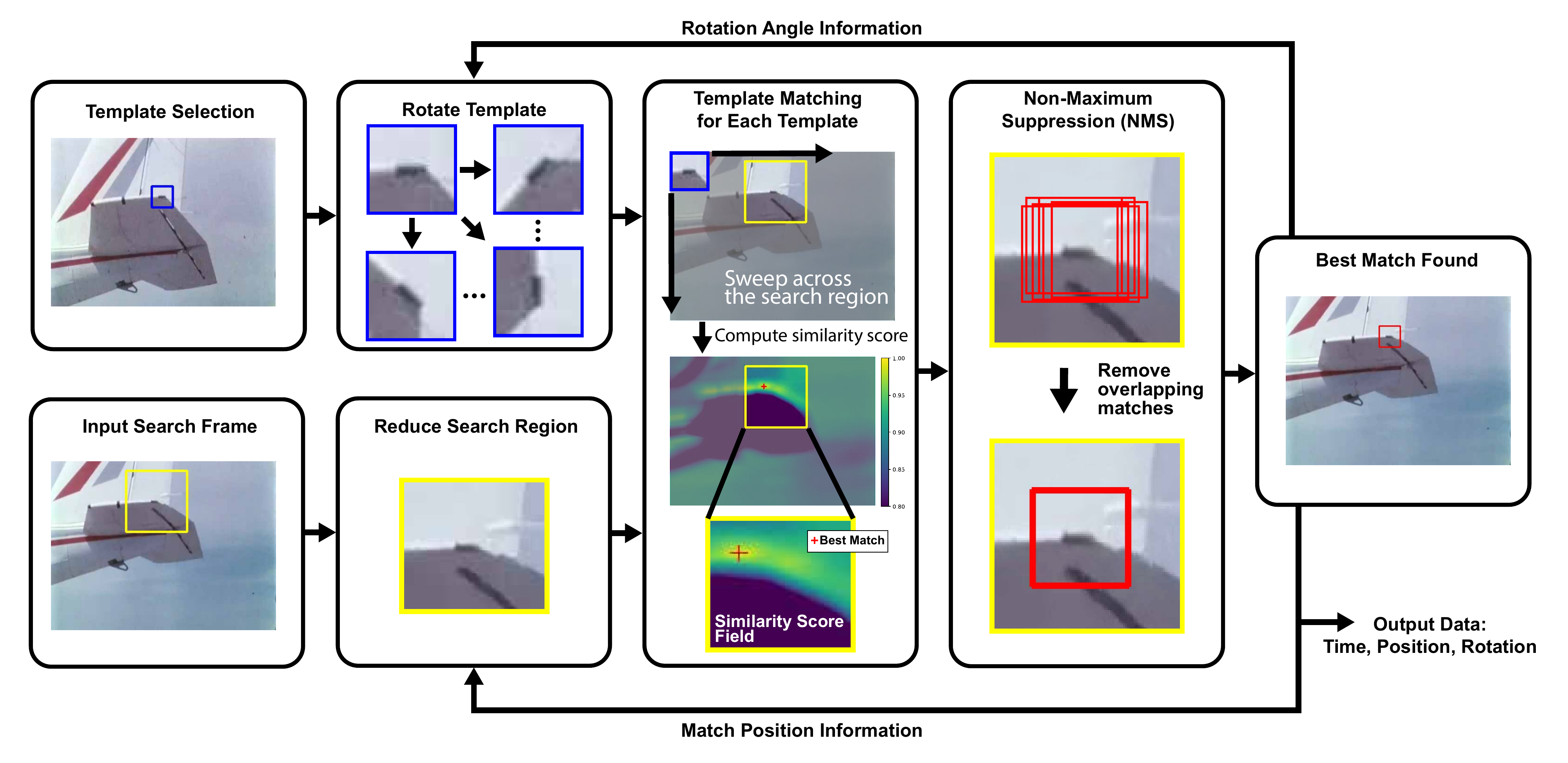}
    \caption[Tracking]%
    {\textbf{Schematic of Template-Matching Tracking Algorithm} --- 
    We initialize one or multiple templates whose centers correspond to the points to track.
    This initialization is only done in the first frame of the video.
    On the next frame of the video, we crop out a smaller search region to look for the template we have selected from the initial frame.
    The template matching algorithm is then applied for each template rotated at fixed intervals, giving optimal match locations for each.
    Note that the similarity score field in the figure shows $1-\bm{R}$ for better visualization.
    The red cross indicates the location of the optimal match. 
    If multiple good matches overlap each other, we remove the redundant matches with the Non-Maximum Suppression algorithm \cite{NMS}.
    With the best match found at the target frame, the match location and rotational angle can be used to update the search region as well as the range of template rotations for the subsequent frames.
    The output from each iteration is collected as time-series data.
    }
    \label{fig:tracking}
\end{figure}

In this section, we describe our model-free tracking algorithm used for recovering trajectory data from experimental videos, which we have implemented in Python using the open-source OpenCV framework \cite{opencv}.
Figure \ref{fig:tracking} illustrates the flowchart of this tracking algorithm which we also summarize later in Algorithm \ref{alg:rotate_tm}.

\subsection{Template matching and tracking algorithm}

Template matching has been actively developed in the computer vision community for the last few decades \cite{tm1_ssd_error,tm2,tm3,tm_driver_fatigue,tm_affine,tm_frame_difference,tm_mean_shift,tm_partition_search,tm_review} with applications to manufacturing \cite{tm_battery,tm_led}, medical imaging \cite{tm_tumor,tm_sperm}, geoinformatics \cite{tm_road}, and general object tracking \cite{tm_real_time,tm_drone,tm_face,tm_rotate_1,tm_rotate_2,tm_underwater_robot}. The technique identifies instances of a predefined template image within a larger target image by comparing pixel-by-pixel similarity \cite{computer_vision}. This is achieved by computing match scores, typically a cross-correlation or sum of square errors, from the intensity values of the template and those of small patches of the target image with the same pixel dimensions. The optimal match is selected as the location with the best match score. This approach is applicable under our assumption that the appearance and size of the object to be modeled remain constant over time. The matching procedure is summarized in lines 5-14 in Algorithm \ref{alg:rotate_tm}.

More specifically, the template matching process starts by selecting a template, a rectangular sub-image, from a reference video frame $\bm{I}_{\textnormal{reference}}(x,y) \in \mathbb{R}^c$, where $c$ is the number of channels.
For a gray-scale image $c=1$ and RGB image $c=3$.
We denote the $i$-th channel of an image with subscript $(\cdot)_i$ (i.e. RGB image $\bm{I} = (I_1, I_2, I_3)$).
The reference video frame is usually the first frame of the video.
Let $\bm{T}(x_t, y_t) \in \mathbb{R}^c$ represent the template image with a separate set of coordinates $(x_t, y_t)$ and $\bm{I}(x,y) \in \mathbb{R}^c$ represent a target frame where we search for the template.
The dimensions, in pixel widths and heights, of the template $\bm{T}$ need to be strictly smaller than the target image $\bm{I}$.
If a background-removed image mask was used, then $T \in \mathbb{B}$ and $I \in \mathbb{B}$, where $\mathbb{B}$ represents the set of Boolean values $\{0,1\}$.

We use image intensity as the matching quantity and the Normalized Sum of Square Difference (NSSD) as the similarity metric \cite{tm1_ssd_error}. By computing a similarity score between the template $\bm{T}$ and all locations on $\bm{I}$ with the NSSD, we obtain a similarity function $R \in \mathbb{R}$ from
\begin{align}
    R(x, y) = r(\bm{T}, \bm{I}) = \frac{\sum_{x_t, y_t, i} \left( T_i(x_t, y_t) - I_i(x_t+x, y_t+y)\right)^2}{\sqrt{\sum_{x_t, y_t, i}T_i(x_t, y_t)^2\sum_{x_t, y_t, i}I_i(x_t+x,y_t+y)^2}},
    \label{eq:NSSD}
\end{align}
where $r$ is an operator that maps a pair of template and target image to a similarity function $R$.
We note that $R(x,y)$ is a scalar function as each summation in (\ref{eq:NSSD}) is done over all $c$ channels.
This process is akin to applying a convolution filter on the target image, where a small template image is swept across the target image and a metric is computed for each location.

The best similarity score (i.e., the smallest value of $R(x,y)$) indicates the optimal match between the template and the search image.
We extend this procedure to account for object rotations by computing similarity functions for templates rotated at fixed intervals and finding the minimum NSSD across all the similarity functions with the Non-Maximum Suppression (NMS) algorithm (described in subsection \ref{section:nms}).
By repeating the matching for all subsequent frames in the video, the algorithm achieves accurate motion tracking over time.

To improve computational efficiency, we utilize the previously detected location and angle to narrow the search region and angle sweep of the next frame.
Thus, instead of processing an entire frame, we search a smaller sub-image $\bm{S}(x,y) \in \mathbb{R}^c$ of $\bm{I}(x,y)$ (see line 4 of Algorithm \ref{alg:rotate_tm}).
By assuming the motion in the video data exhibits spatial and temporal coherence, we limit the search area to a window surrounding the previously matched location and rotation range.
The size of the search window and the range of rotations must exceed the maximum displacement and rotation observed in the tracked object throughout the video.
Therefore, having prior knowledge of the system and video quality (such as frame rate) can inform the selection of these parameters.
Based on the examples we presented in this work, a search region of twice the pixel length of the template dimensions and a rotation range within $\pm 15 ^{\circ}$ degrees about previously matched template rotation with an interval of $5 ^{\circ}$ is a good starting point.

\subsection{Background removal through frame averaging}\label{section:bg_subtract}
A helpful preprocessing step for motion extraction involves separating foreground objects from the background.
A simple yet effective technique is background subtraction through frame averaging.
This preprocessing step is applied to the double pendulum and inverted flag examples presented in Section \ref{section:examples}.
The technique corresponds to lines 1-4 in Algorithm \ref{alg:rotate_tm}.

Frame averaging seeks to extract the foreground objects from a sequence of images or video by subtracting the static background, obtained from time-averaging image intensities from each video frame.
This method works under the assumption that the background in a video sequence tends to remain constant or have minimal changes over time, while the foreground objects introduce significant variations in terms of image intensities.

The process of frame averaging for an RGB-colored video $\bm I\in \mathbb{R}^3$ is represented by
\begin{align}
    \bm{I}_\mathrm{mean}(\bm{x}) = \frac{1}{N} \sum_{t} \bm{I}(\bm{x}, t),
\end{align}
where $\bm{x} = (x,y)$, $t$ is the time corresponding to each frame, and $N$ is the total number of frames in the video.
This operation is allowed as the RGB color space is a linear additive space.

The average frame $\bm{I}_\mathrm{mean}$ is subtracted from each individual frame in the video sequence.
The result is a difference image $\bm{D}\in\mathbb{R}^3$,
\begin{align}
    \bm{D}(\bm{x}, t) = |\bm{I}(\bm{x}, t) - \bm{I}_\mathrm{mean}(\bm{x})|,
\end{align}
that highlights the foreground objects.
The larger the magnitude of the difference, the more likely the pixel is in the foreground.
 
To refine the foreground-background segmentation, we take the maximum value of the $c$ channels of the difference frame, that is
\begin{align}
    D_\mathrm{max}(\bm{x}, t) =  \max_{1\leq i\leq c}{D_i(\bm{x}, t)}.
\end{align}

By thresholding $D_\mathrm{max}$ as
\begin{align}
    D_\mathrm{mask}(\bm{x}, t) = h(D_\mathrm{max}, D_\mathrm{thresh}) = \begin{cases}
                                    1, &D_\mathrm{thresh} \leq D_\mathrm{max}(\bm{x}, t)\\
                                    0, &\mathrm{otherwise}
                              \end{cases},
\end{align}
we obtain a boolean mask that separates the foreground from the background.
The mask contains foreground shape information and we can track features directly in this binary image.

\subsection{Non-maximum suppression}\label{section:nms}
To remove redundant bounding boxes or detections, we use the non-maximum suppression (NMS) algorithm \cite{NMS}.
It involves sorting the bounding boxes based on their confidence scores, selecting the box with the highest score, and suppressing other boxes with high overlap.
This process ensures that only the most accurate and non-overlapping detections are retained.
NMS is widely employed in object detection algorithms to improve accuracy by eliminating duplicate detections.
As it is expected that our problems have multiple overlapping close matches from different rotated templates, we use the NMS algorithm to eliminate the ambiguity systematically.
In all the examples presented in this work, we keep only the global best match.

The overlap of the bounding boxes is measured with the Jaccard index \cite{Jaccard}, which is more commonly known as the Intersection over Union (IoU) metric in computer vision, defined as
\begin{align}
    J(A,B) = \left| \frac{A \cap B}{A \cup B}\right| = \frac{\left|A \cap B\right|}{|A| + |B| - \left|A \cap B\right|},
    \label{eq:Jaccard}
\end{align}
where $|\cdot|$ indicates the size of set, and $J(A,B)$ has a range of $[0,1]$ by design.

The NMS algorithm corresponds to lines 15-30 in Algorithm \ref{alg:rotate_tm} and is summarized in the following steps:
\begin{enumerate}
    \item Obtain a list of detection boxes $\mathcal{L}$ with a corresponding list of confidence scores $\mathcal{R}$ and rotations $\mathcal{T}$.
    \item Find the detection, $\mathcal{L}_m$ with the best similarity score (minimum of $\mathcal{R}$) from set $\mathcal{L}$.
    \item Remove the detection $\mathcal{L}_m$ from the set $\mathcal{L}$ and appends it to the set of final output detections $\mathcal{K}$.
    Repeat this for match scores $\mathcal{R}$ and rotations $\mathcal{T}$.
    \item Using the Jaccard index (\ref{eq:Jaccard}) as the overlap metric, remove all detection boxes that have an overlap with $\mathcal{L}_m$ greater than a threshold $\mathrm{IoU}_\mathrm{thresh}$ in the set $\mathcal{L}$.
    Repeat this for match scores $\mathcal{R}$ and rotations $\mathcal{T}$.
    \item Repeat step 2-4 until either $\mathcal{L}$ is empty or the size of $\mathcal{K}$ reaches the required number of non-intersecting matches $n_\mathrm{matches}$.
\end{enumerate}
As we track a single point in all five video examples, we set $n_\mathrm{matches} = 1$.
We remark that our proposed tracking algorithm with NMS can potentially track multiple identical objects, however this is not pursued in this paper.
Lastly, we find an overlap threshold of $\mathrm{IoU}_\mathrm{thresh} = 0.3$ to be effective across all videos.

\subsection{Summary}
Our template-matching-based tracking algorithm is summarized in Algorithm \ref{alg:rotate_tm}.
The tracking algorithm has three main steps: 1. Optional background subtraction pre-processing, 2. Generate candidate templates and apply template matching, 3. Apply non-maximum suppression to obtain optimal estimate of object location.
We apply the algorithm to every video frame and collect the output position and rotation of the tracked template as time-series data.

The inputs to the algorithm/tracking routine are:
\begin{enumerate}
    \item user-selected template in the form of a small patch on the first frame of a video,
    \item target frame to search for matches,
    \item optional averaged frame to use for background subtraction,
    \item threshold value for maximum difference image to refine background removal, 
    \item a region on the target frame to search for matches which is updated at every frame based on the previous match location,
    \item intersection over union threshold to remove the overlapping bounding box from matches,
    \item angle sweep bounds and intervals which are updated at every frame based on the rotation previous template match,
    \item number of non-overlapping best matches to search for.
\end{enumerate}

In practice, the user only needs to initialize the tracker by selecting a region, to be tracked, in the first frame.
The tracking algorithm is evaluated at every frame of the video, outputting the $(x,y)$ locations of the best matches as well as the corresponding rotation angle of the template. 
The algorithm then uses the match location information to update the search region and the angle sweep bounds for the next frame to speed up computation.
We find updating the search region to a region within 2 times the pixel length of the template centering the previously matched $(x,y)$ location and $\pm 15 ^{\circ}$ about the previously matched rotation angle with an interval of $5 ^{\circ}$ work well for the examples presented in this work.

\renewcommand{\algorithmicrequire}{\textbf{Input:}}
\renewcommand{\algorithmicensure}{\textbf{Output:}}
\newlength{\commentlen}
\setlength{\commentlen}{40ex}
\newcommand{\algComment}[1]{\Comment{\makebox[\commentlen][l]{#1}}}

\begin{algorithm}[H]
\caption{Template Matching-Based Tracking Algorithm}\label{alg:rotate_tm}
\begin{algorithmic}[1]

\Require
    \Statex $\bm{T}(x_t,y_t)$, template image
    \Statex $\bm{I}(x,y)$, target frame
    \Statex $\bm{I}_\mathrm{mean}(x,y)$, background image from frame average
    \Statex $D_\mathrm{thresh}$, difference image threshold value
    \Statex \textit{search\_region}, region on target frame for matching
    \Statex $\mathrm{IoU}_\mathrm{thresh}$, Intersection over Union (overlap) threshold
    \Statex $\theta_\mathrm{min}$, $\theta_\mathrm{max}$, $\theta_\mathrm{interval}$, angle sweep bounds and interval
    \Statex $n_\mathrm{match}$, number of best matches
\Ensure
    \Statex $\mathcal{K}$, list of top left $(x,y)$ coordinates of found matches
    \Statex $\mathcal{T}$, list of rotation of the found matches

\Statex

\Statex $\text{/* preprocess by removing background */}$
    \State $\bm{D} \gets |\bm{I} - \bm{I}_\mathrm{mean}|$
    \State $D_\mathrm{max} \gets \max_{1\leq i\leq3}{D_i}$
    \State $D_\mathrm{mask} \gets h(D_\mathrm{max}, D_\mathrm{thresh})$ \algComment{threshold}
    \State $\bm{S}\gets$ crop $\bm{I}$ or $D_\mathrm{mask}$ at \textit{search\_region}

\Statex $\text{/* template matching */}$
    \State $\mathcal{L} \gets \{\}$ \algComment{match position array}
    \State $\mathcal{R} \gets \{\}$ \algComment{match score array}
    \State $\mathcal{T} \gets \{\}$ \algComment{match angle array}

    \For{$\theta \in \{x: \theta_\mathrm{min}\leq x \leq \theta_\mathrm{max}, \, \frac{x-\theta_\mathrm{min}}{\theta_\mathrm{interval}}\in \mathbb{Z}\}$}
        \State rotate $\bm{T}$ by $\theta$
        \State $\bm{R} \gets r(\bm{T}, \bm{S})$ \algComment{template match}
        \State $\mathcal{L} \gets \mathcal{L} \cup \{(x,y) \mid \bm{R}(x,y) < 0.5\}$
        \State $\mathcal{R} \gets \mathcal{R} \cup \{(x,y) \mid \bm{R}(x,y) < 0.5\}$
        \State $n \gets |\{(x,y) \mid \bm{R}(x,y) < 0.5\}|$ \algComment{number of potential good matches}
        \State $\mathcal{T} \gets \mathcal{T} \cup \{\theta\}_{n}$
    
\algstore{myalg} 
\end{algorithmic}
\end{algorithm}

\begin{algorithm}[H]
\begin{algorithmic}[1]
\algrestore{myalg} 
    \Statex $\text{/* Non Maximum Suppression */}$
        \State $\mathcal{K} \gets \{\}$ \algComment{match position array for output}
        \While{$\mathcal{L} \neq \emptyset \textnormal{ and } |\mathcal{K}| < n_\textnormal{match}$}
            \State $m\gets \mathrm{argmin}_j \mathcal{R}_j$
            \State $\mathcal{K} \gets \mathcal{K} \cup \mathcal{L}_m$
            \State $\mathcal{L} \gets \mathcal{L} - \mathcal{L}_m$
            \State $\mathcal{R} \gets \mathcal{R} - \mathcal{R}_m$
            \State $\mathcal{T} \gets \mathcal{T} - \mathcal{T}_m$
    
            \For{$\mathcal{L}_i \; \in \; \mathcal{L}$}\algComment{remove overlapping matches}
                \If{$J(\mathcal{L}_m, \mathcal{L}_i) \geq \mathrm{IoU}_\mathrm{thresh}$}
                    \State $\mathcal{L} \gets \mathcal{L} - \mathcal{L}_i$
                    \State $\mathcal{R} \gets \mathcal{R} - \mathcal{R}_i$
                    \State $\mathcal{T} \gets \mathcal{T} - \mathcal{T}_i$
                \EndIf
            \EndFor
        \EndWhile
\EndFor

\State \textbf{return} $\mathcal{K}$, $\mathcal{T}$
    
\end{algorithmic}
\end{algorithm}

\section{Data-driven reduced-order models on spectral submanifolds}
Here we summarize available results on rigorous model order reduction to SSMs in smooth nonlinear systems and the \emph{SSMLearn} algorithm we have used for constructing SSM-reduced models from data.

\subsection{Problem setup}
Consider an $n$-dimensional dynamical system of the form
\begin{align}\label{eq:system}
    \dot{\bm{x}} &= \bm{Ax} + \bm{f}(\bm{x}),\\
    \bm{x}&\in \mathbb{R}^n, \; \bm{A}\in\mathbb{R}^{n\times n}, \; \bm{f}:\mathbb{R}^n \mapsto\mathbb{R}^n, \; \bm{f} \sim \mathcal{O}(|\bm{x}|^2), \; \bm{f}(\bm{0}) = \bm{0}, \nonumber
\end{align}
where $\bm{f}$ is of differentiability class $C^r$ in $\bm{x}$ and $\bm{A}$ is the linear part of the system.

For simplicity of exposition, we assume that $\bm{A}$ is diagonalizable and $\bm{x} = \bm{0}$ is an asymptotically stable fixed point, that is
\begin{equation}
    \mathrm{Re}\,\lambda_k<0\quad \forall \lambda_k\in\mathrm{Spect}(\bm{A}),
\end{equation}
where $\mathrm{Spect}(\bm{A}) = \{\lambda_1, \lambda_2, ..., \lambda_n \}$ denotes the eigenvalues of the linear part of the system at the origin.
Recent work extended this result from stable to general hyperbolic fixed points with possible mix of stable and unstable eigenvalues \cite{ssm_mixed}.

We select a $d$-dimensional spectral subspace $E \subset \mathbb{R}^n$, that is, the direct sum of a set of eigenspaces of $\bm{A}$.
We denote the spectrum of $\bm{A}$ within $E$ by $\mathrm{Spect}(\bm{A}|_E)$.
If the spectral subspace $E$ is non-resonant with the spectrum of $\bm{A}$ outside of $E$, i.e.,
\begin{align}
    \lambda_j &\neq \sum^{d}_{k=1} m_k\lambda_k, \quad
    m_k \in \mathbb{N}, \quad
    \sum^{d}_{k=1} m_k \geq 2,
    \\
    \lambda_j &\in \mathrm{Spect}(\bm{A}) \backslash \mathrm{Spect}(\bm{A}|_E), \quad
    \lambda_k \in \mathrm{Spect}(\bm{A}|_E),
    \nonumber
\end{align}
then in general $E$ has infinitely many nonlinear continuations, as invariant manifolds of dimension $d$ emanating from $\bm{x}=\bm{0}$ and tangent to $E$, in the system (\ref{eq:system}) \cite{ssm}.
Within this family of invariant manifolds, there is a unique smoothest member, which we define as the primary spectral submanifold.
This SSM is normally attracting, can be computed as a Taylor expansion, and is therefore an ideal candidate for nonlinear model reduction.
SSMs can be efficiently computed from the equations of motion using the open-source package \emph{SSMTool} \cite{SSMTool2}.

\subsection{SSM-based data-driven modeling algorithm}
When the equations of motion are unknown, model reduction can also be applied directly to time-dependent observable data describing the evolution of the dynamical system.
Recently, SSM theory has been applied to both simulated and experimental data to capture essential nonlinear dynamics and enable accurate predictions for motions not used in the training or for behavior under additional external forcing \cite{ssm_data_old,SSMLearn,ssm_soft_robot,liu24}.
Such data-driven models can even outperform analytical models as they are not necessarily bound to the domain of convergence of the Taylor expansions used for the SSMs and their  reduced dynamics  \cite{fastSSM}.

To construct SSM-reduced models from data, we use the methodology presented in \cite{SSMLearn}, which is implemented in the open-source MATLAB package \emph{SSMLearn}.
The method consists of two main steps: geometry identification and reduced dynamics modeling.
As the low-dimensional SSM attracts nearby trajectories, we approximate it as a polynomial using nearby transient training data in the phase space or observable space.
We then model the reduced dynamics by identifying a vector field or map from the high-dimensional training trajectories projected onto the SSM local coordinates.
Here we provide a summary of the \emph{SSMLearn} algorithm; a complete description can be found in \cite{SSMLearn}.

We start with the geometry identification.
To obtain a graph-style model of a $d$-dimensional SSM from data, we construct the manifold as a polynomial parameterized by local coordinates on its tangent space. 
To this end, we define a matrix $\bm{V} \in \mathbb{R}^{n\times d}$, whose columns are orthonormal vectors spanning the tangent space of the yet unknown SSM.
The reduced coordinates $\bm{\xi}(t) \in \mathbb{R}^d$ are defined as a projection of the trajectories $\bm{y}(t) \in \mathbb{R}^n$ onto the tangent space $\bm{V}$, that is
\begin{equation}\label{eq:ssm-tangent-space}
    \bm{\xi} = \bm{V}^{\top} \bm{y}.
\end{equation}

We seek to approximate and parameterize the manifold with a Taylor expansion in $\bm{\xi}$ about the fixed point at the origin, $\bm{y} = \bm{0}$
\begin{equation}\label{eq:ssm-manifold-parameterisation}
    \bm{y} \approx \bm{M} \bm{\xi}^{1:m} = \bm{V}\bm{\xi} + \bm{M}_{2:m}\bm{\xi}^{2:m} = \bm{v}(\bm{\xi}),
\end{equation}
where $\bm{M}\in \mathbb{R}^{n\times d_{1:m}}$ is a matrix of manifold parameterization coefficients, with $d_{1:m}$ denoting the number of $d$-variate monomials from orders 1 up to $m$. $(\cdot)^{l:m}$ denotes a vector of all monomials at orders $l$ through $m$. For example, if $\bm{x} = [x_1, x_2]^{\top}$, then $\bm{x}^{2:3} = [x_1^2, x_1 x_2, x_2^2, x_1^3, x_1^2 x_2, x_1 x_2^2, x_2^3]^{\top}$.
$\bm{M}_{2:m}$ denotes a submatrix of $\bm{M}$ that has only columns associated with $(\cdot)^{2:m}$ vector, so $\bm{M}_{2:m} \in \mathbb{R}^{n\times d_{2:m}}$.

The matrices $\bm{V}$ and $\bm{M}$ are solved for simultaneously by minimizing the cost function
\begin{equation}
    (\bm{V}, \, \bm{M}) = \operatorname*{argmin}_{(\bm{V}^*, \, \bm{M}^*)} \sum_{j=1}^N \|\bm{y}_j - \bm{M}^*(\underbrace{\bm{V}^{*\top} \bm{y}_j}_{\bm{\xi^*}})^{1:m} \|^2,
\end{equation}
subject to constraints
\begin{equation}
    \bm{V}^{\top}\bm{V} = \bm{I}, \; \bm{V}^{\top}\bm{M}_{2:m} = \bm{0}.
\end{equation}

In the second step, after we have found the parameterization of the manifold, we approximate the reduced dynamics up to $r$-th order on the manifold as
\begin{equation}
    \dot{\bm{\xi}} \approx \bm{R} \bm{\xi}^{1:r} = \bm{r}(\bm{\xi}),\\
\end{equation}
where $\bm{R}\in\mathbb{R}^{d\times d_{1:r}}$ is a matrix of coefficients of the reduced dynamics.
We solve for the matrix $\bm{R}$ through minimization of the cost function
\begin{equation}
    \bm{R} = \operatorname*{argmin}_{\bm{R}^*} \sum_{j=1}^N \| \dot{\bm{\xi}}_j - \bm{R}^* \bm{\xi}_i^{1:r} \|^2 = \operatorname*{argmin}_{\bm{R}^*} \sum_{j=1}^N \left|\left| \frac{d}{dt}(\bm{V}^{\top}\bm{y}_j) - \bm{R}^* (\bm{V}^{\top}\bm{y}_j)^{1:r} \right|\right|^2.
\end{equation}

We then transform the reduced dynamics to its normal form, which is the simplest complex polynomial form that preserves the differential topology of the local dynamics.
Transforming the equations to a normal form therefore provides insights into the qualitative behaviors of the system, such as stability properties, fixed points, limit cycles, and bifurcations \cite{Guckenheimer}.

We compute the $n$-th order normal form of $\bm{R}$ to obtain a near-identity polynomial transformation of $\bm{\xi}$ to complex conjugate normal form coordinates $\bm{z} \in \mathbb{C}^{d}$.
We define the transformations $\bm{t}:\bm{z}\mapsto \bm{\xi}$, its inverse $\bm{t}^{-1}:\bm{\xi}\mapsto \bm{z}$, and the dynamics in normal form $\bm{n}:\bm{z}\mapsto \dot{\bm{z}}$ as
\begin{alignat}{4}
    \bm{\xi} &= \bm{t}(\bm{z}, \bm{T}) = \bm{T}\bm{z}^{1:n} = \bm{W}\bm{z} + \bm{T}_{2:n}\bm{z}^{2:n},\label{eq:ssm-z_to_xi} \\
    \bm{z} &= \bm{t}^{-1}(\bm{\xi}, \bm{H}) = \bm{H}\bm{\xi}^{1:n} = \bm{W}^{-1}\bm{\xi} + \bm{H}_{2:n}(\bm{W}^{-1}\bm{\xi})^{2:n}, \label{eq:ssm-xi_to_z}\\
    \dot{\bm{z}} &= \bm{n}(\bm{z}, \bm{N}) = \bm{N}\bm{z}^{1:n} = \bm{\Lambda z} + \bm{N}_{2:n}\bm{z}^{2:n},\label{eq:ssm-z_dynamics}
\end{alignat}
\begin{align*}
    \bm{T} \in \mathbb{C}^{d\times d_{1:n}}, \; \bm{H} \in \mathbb{C}^{d\times d_{1:n}}, \; \bm{N} \in \mathbb{C}^{d\times d_{1:n}},
\end{align*}
where $\bm{\Lambda}\in\mathbb{R}^{d\times d}$ is a diagonal matrix of eigenvalues of the linear part $\bm{R}_{1:1} = \bm{W\Lambda W}^{-1}$ of $\bm{R}$, and $\bm{W}\in\mathbb{R}^{d\times d}$ contains the associated eigenvectors.

The normal form transformation and its dynamics are solved for simultaneously with the minimization of the cost function
\begin{equation}
    (\bm{H}, \, \bm{N}) = \operatorname*{argmin}_{(\bm{H}^*, \, \bm{N}^*)} \sum_{j=1}^N \| \underbrace{\nabla_{\bm{\xi}} \bm{t}^{-1}(\bm{\xi}_j, \bm{H}^*)\dot{\bm{\xi}}_j}_{\dot{\bm{z}}^*} - \bm{n}(\underbrace{\bm{t}^{-1}(\bm{\xi}_j, \bm{H}^*)}_{\bm{z}^*}, \bm{N}^*) \|^2.
\end{equation}

Finally, with the the matrix $\bm{H}$ computed, we solve for $\bm{T}$ with 
\begin{equation}
    \bm{T} = \operatorname*{argmin}_{\bm{T}^*} \sum_{j=1}^N \| \bm{t}(\underbrace{\bm{t}^{-1}(\bm{\xi_j}, \bm{H})}_{\bm{z}_j}, \bm{T^*}) - \bm{\xi}_j \|^2.
\end{equation}

Our reduced model is thus described fully by (\ref{eq:ssm-tangent-space})(\ref{eq:ssm-manifold-parameterisation})(\ref{eq:ssm-z_to_xi})(\ref{eq:ssm-xi_to_z})(\ref{eq:ssm-z_dynamics}).

\subsection{Delay embedding}\label{section:delay_embed}
In practice, observing or measuring all variables that span the full phase space of a dynamical systems is impractical if not impossible. This prompts us to construct  data-driven SSM-reduced models in an appropriate observable space. We achieve this by delay-embedding the SSM in the observable space based on the Takens embedding theorem \cite{Takens}.

Delay embedding involves constructing an embedding space using time-delayed copies of a given observable scalar time series $s(t) = \mu(\bm{x})$, where $\mu$ is a differentiable scalar observable function $\mu: \mathbb{R}^{n_\mathrm{full}} \mapsto \mathbb{R}$.
The function $\mu$ returns a measured feature of the system (\ref{eq:system}), such as the displacement of a particular material point.
A new state vector $\bm{y} \in \mathbb{R}^p$ is constructed from $p$ delayed measurements of the original scalar time series separated by a timelag $\Delta t$, that is
\begin{align}
    \bm{y}(t) = \begin{bmatrix}
                s(t) \\
                s(t+\Delta t) \\
                s(t+2\Delta t) \\
                \vdots \\
                s(t+(p-1)\Delta t) 
            \end{bmatrix}.
\end{align}

Takens's theorem states that if $\bm{x}(t)$ lies on an $d$-dimensional invariant manifold in the full phase space, then the embedded state $\bm{y}(t)$ also lies on a diffeomorphic manifold in the $p$-dimensional embedding space, if $p\geq 2d + 1$ and $\mu$ and $\Delta t$ are generic (as defined in Remark 1 in \cite{ssm_delay_embed_review}).
This result can also be extended to higher-dimensional observable quantities ($\bm{\mu}: \mathbb{R}^{n_\mathrm{full}} \mapsto \mathbb{R}^{n_\mathrm{observe}}$), under appropriate nondegeneracy conditions outlined in \cite{ssm_delay_embed_review, delay_embed_1, delay_embed_2}.

Near a fixed point, more can be said about the structure of the embedding.
It was recently shown that the orientation of eigenspaces at a fixed point in delay-embedded space is directly determined by the eigenvalues of the system linearized there \cite{ssm_delay_embed_review}.
Thus local spectral properties of the full phase space have a direct geometrical interpretation in the observable space that can be exploited to improve system identification.

In our context, we interpret a video as a generic observable function of the dynamical system.
The coordinates we extract from it are in turn observables on that video.
The extracted pixel data is therefore an observable of the full system and Takens's theorem applies after an appropriate delay-embedding.

\subsection{Backbone curve from the normal form}\label{section:backbone}

By expressing $\bm{z}$ in polar coordinates $(\bm{\rho}, \bm{\theta})$, defined as $z_j = \rho_j e^{i\theta_j}$ for $j=1,...,m$, we obtain the extended normal form of an oscillatory SSM of dimension $2m$ of the form
\begin{align}
    \dot{\rho}_j &= \gamma_j(\bm{\rho}, \bm{\theta})\rho_j,\label{eq:general_backbone_rho} \\
    \dot{\theta}_j &= \omega_j(\bm{\rho}, \bm{\theta}),\label{eq:general_backbone_theta} \\
    j &= 1,...,m, \; \bm{\rho}\in\mathbb{R}^m_+, \; \bm{\theta}\in \mathbb{T}^m, \nonumber
\end{align}
which enables us to perform analyses such as distinguishing different modal contributions, applying slow-fast decomposition, and constructing backbone curves.

If there is no resonance in the linearized eigenfrequencies, then $\gamma_j$ and $\omega_j$ from (\ref{eq:general_backbone_rho}) and (\ref{eq:general_backbone_theta}) are functions of the amplitudes $\bm{\rho}$ only \cite{SSMTool1}.
This allows us to decouple the amplitude dynamics from the phase dynamics.

The zero-amplitude limit of $\gamma_j$ and $\omega_j$ corresponds to the linearized damping and frequency of the $j$th eigenmode, that is
\begin{align}
    \lim_{||\bm{\rho}|| \rightarrow \bm{0}} (\gamma_j(\bm{\rho}, \bm{\theta}) + i\omega_j(\bm{\rho}, \bm{\theta})) = \lambda_j.
\end{align}
Thus, $\gamma_j$ and $\omega_j$ represent the nonlinear extension of damping and frequency of the system.

To gain physical insight into the amplitude dependence of instantaneous damping and frequency through backbone curves, the amplitude in the normal coordinates must be transformed back into the observable space.
For a general 2D case, the amplitude $\mathcal{A}$ can be determined through 
\begin{align}
    \mathcal{A}(\rho) = \max_{\theta \in [0,2\pi)} |g(\bm{v}(\bm{t}(\bm{z})))|, \; \bm{z}=\left( \rho e^{i\theta}, \rho e^{-i\theta}\right),
\end{align}
where the function $g:\mathbb{R}^n \mapsto \mathbb{R}$ maps from the observable space to the amplitude of a particular observable, such as a degree of freedom.

The backbone and damping curves are expressed as parametric functions
\begin{align}
    \mathcal{B}_{\mathrm{damping}} &= \{
                                        \gamma(\rho),
                                        \mathcal{A}(\rho)
                                         \}_{\rho\geq 0}, \label{eq:backbone_damp}\\
    \mathcal{B}_{\mathrm{frequency}}& =  \{\omega(\rho),\mathcal{A}(\rho)\}_{\rho\geq 0}.
                                            \label{eq:backbone_frequency}
\end{align}

\section{Applications}\label{section:examples}
We now apply our tracking and \emph{SSMLearn} algorithms to five examples with little to no modifications to the video tracking and model training procedure.

\subsection{Double pendulum}
\begin{figure}[t]
    \centering
    \includegraphics[width=1\linewidth]{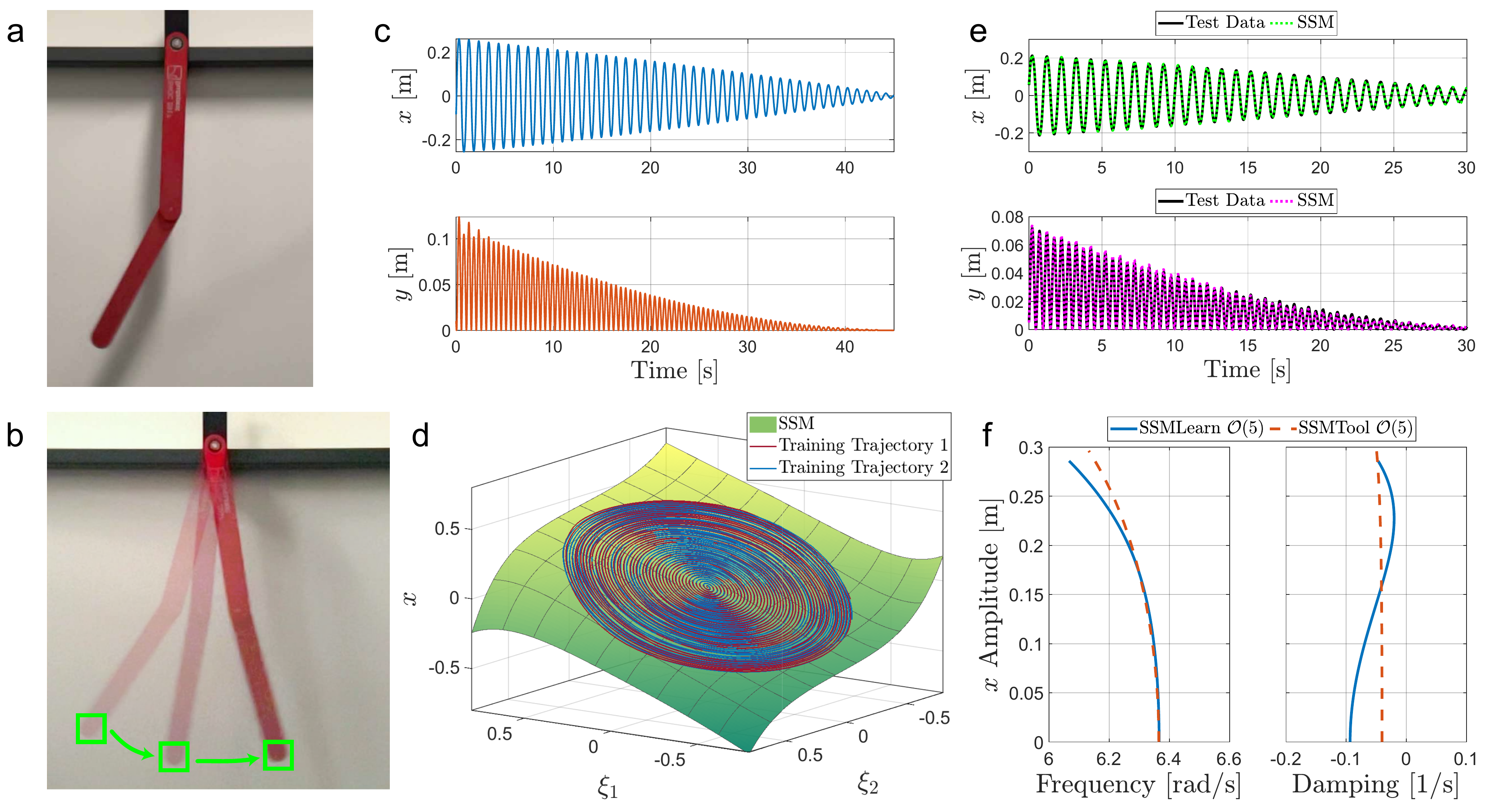}
    \caption[DP]%
    {\textbf{Data-driven Nonlinear Reduced-order Model on the Slowest SSM of a Double Pendulum} --- 
    (a) System setup of the double pendulum.
    (b) Snapshots of processed video frames.
    (c) Training trajectory in physical length, which was computed by scaling pixels to known physical length from measurements.
    (d) Training trajectory and fitted manifold plotted in the two reduced coordinates $\xi_{1,2}$ and $x$.
    (e) Test trajectory and its reconstruction by the SSM model.
    (f) Backbone (left) and damping (right) curve output of the \emph{SSMLearn} model, showing the nonlinearities in the system. The analytical result computed from \emph{SSMTool} \cite{SSMTool2} is plotted together for reference.
    }
    \label{fig:double_pendlum}
\end{figure}

As our first example, we consider a compound double pendulum as shown in Figure \ref{fig:double_pendlum}(a).
The double pendulum in our study consists of a pair of rods; the upper rod weighs 253 g and measures 200 mm while the lower rod is 180 mm long and weighs 114 g.
The rods are connected with a roller bearing and both have a width of 25 mm.

We film videos of the transient decay response of the double pendulum released from different initial conditions.
All videos are recorded with the iPhone 12 Pro, utilizing only the wide camera, at a frame rate of 240 FPS with 1920x1080 resolution.
The camera has up to 12 megapixels resolution with f/1.6 aperture.

To extract data from the videos, we initialize the tracker once by selecting a template, in the form of a bounding box, in the first frame of the video.
The center of the template corresponds to the end of the lower pendulum rod.
The algorithm tracks the template for all subsequent frames in the video and outputs time-series data of the $(x,y)$ pixel locations of the template center.
Figure \ref{fig:double_pendlum}(b) illustrates snapshots of tracked video.
One of the two training trajectories is shown in panel (c) of Figure \ref{fig:double_pendlum}, while the black trajectories in panel (e) are used as test data. 
Since the physical lengths of the pendulum arms are known, we scale the output pixel data to physical length in meters.

From the video-extracted trajectory data, we aim to learn the reduced dynamics on the slow 2-dimensional SSM emanating from the slow eigenspace of the linear part of the system.
The minimal embedding dimension is 5 according to Takens's theorem ($p\geq 2d+1$).
Thus we embed the training data (both $x$ and $y$) five times with a time step of 0.004167s (240 FPS) to reconstruct the SSM in the observable space.

Figure \ref{fig:double_pendlum}(d) shows the identified 5th-order, 2-dimensional SSM in three coordinates: $x$ and the two reduced coordinates $\xi_{1,2}$.
The choice of polynomial order of the SSM parameterization is a compromise between model accuracy, complexity, and the risk of over-fitting.
For all the examples in the present paper, we choose a manifold expansion order that has less than 3\% parameterization error (see Appendix for error quantification).

On the SSM, \emph{SSMLearn} returns reduced dynamics in a 5th-order normal form
\begin{equation}
    \begin{aligned}
        \dot{\rho}\rho^{-1} &= -0.09352 + 0.8130\rho^2 - 2.256\rho^4,\\
        \dot{\theta} &= 6.366 - 0.4733\rho^2 - 1.953\rho^4.
    \end{aligned}\label{eq:tm_normal_form}
\end{equation}

We now use the model (\ref{eq:tm_normal_form}), trained on a single trajectory, to make predictions of the test trajectory released from a different initial condition.
The reduced model takes the initial condition of the test trajectory as the only input and reconstructs the entire decay response with a reconstruction error of 2.1\%.

For comparison, we also derived the equations of motions for this double pendulum system  and applied the equation-driven SSM identification package,  \emph{SSMTool}, (see the Appendix for more details).
In this calculation, the damping from the two joints was modeled as linear using the Rayleigh dissipation function. The geometry of the double pendulum was modeled by two rods with distributed mass. We constrained the motion of the double pendulum to a vertical plane which resulted in a two-degree-of-freedom mechanical model.
This four-dimensional dynamical system is fully described by two angles $\theta_i$ for $i = 1, 2$ relative to the vertical, and their respective angular velocities $\dot{\theta_i}$.

The frequency and backbone curves obtained from the video data and from the equation of motions are plotted together in Fig. \ref{fig:double_pendlum}(f).
The frequency learned from data is consistent with the analytical result within 2\% error.
However, the analytically predicted nonlinear damping curve on the right differs from the actual damping curve identified from the video. This shows that the idealized linear damping  used in our mechanical model fails to capture the actual low-amplitude damping characteristic of the system, which is likely dominated by dry friction  that is independent of rotational speed.
The presence of dry friction is also evident from the more pronounced model mismatch at lower oscillation amplitude.
This result illustrates well the use of our video-based SSM-reduction procedure in accurate nonlinear system identification.

\subsection{Inverted flag}
\begin{figure}[t]
    \centering
    \includegraphics[width=1\linewidth]{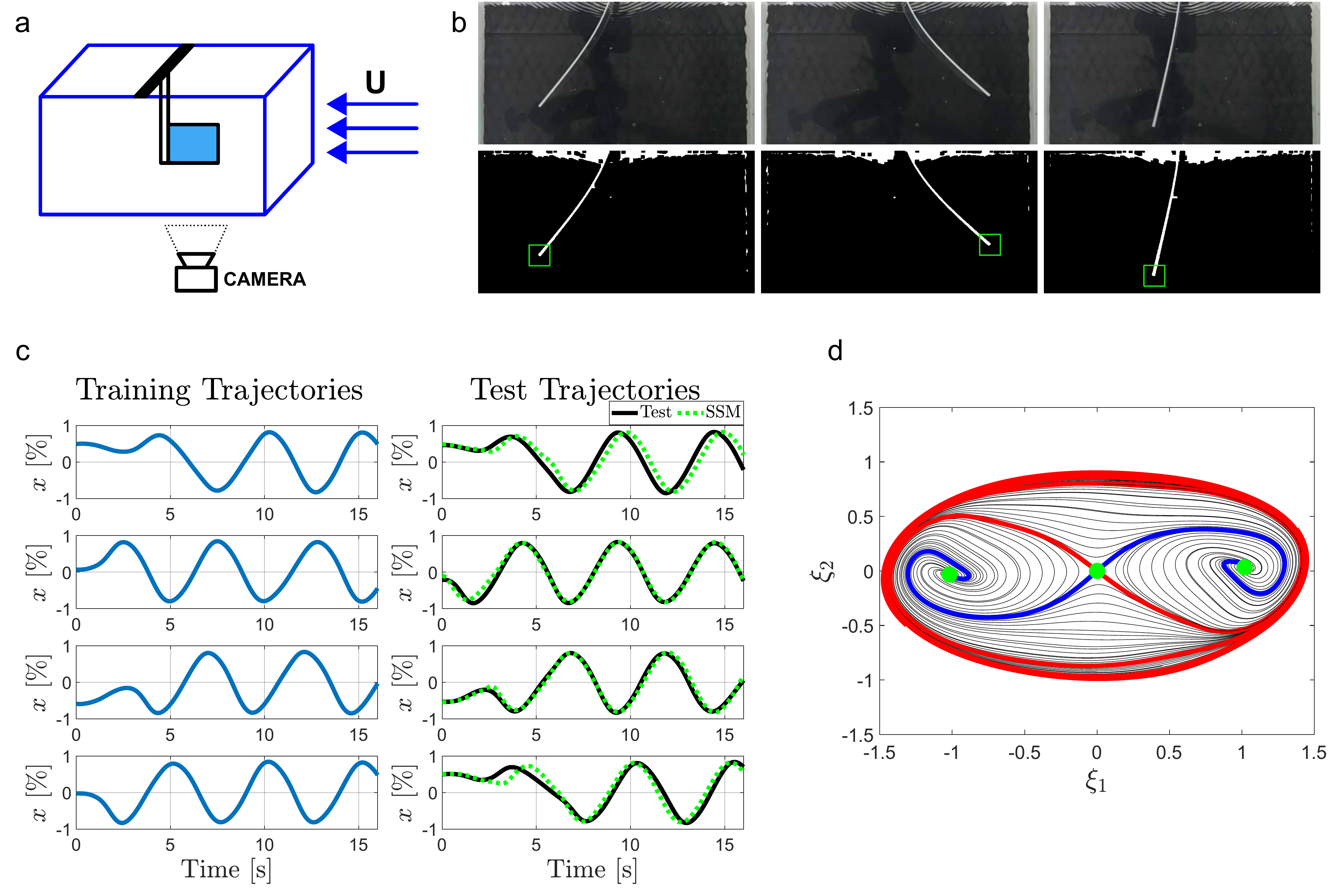}
    \caption[Flag]%
    {\textbf{Data-driven Nonlinear Reduced-order Model on the Mixed-Mode SSM of an Inverted Flag} --- 
    (a) System setup.
    (b) Three selected video frames (top row) and the corresponding background-removed and tracked location (bottom row).
    (c) The training trajectories and test trajectories with model predictions.
    (d) Phase portrait in reduced coordinates, illustrating the three fixed points in green, stable manifold in blue, and unstable manifold in red.
    }
    \label{fig:flag}
\end{figure}

For our second example, we consider an inverted flag, a flexible elastic sheet with its trailing edge clamped subject to a fluid flow.
This is in contrast to a conventional flag where the leading edge is constrained.
The experimental configuration, based on \cite{xu24}, is shown in \ref{fig:flag}(a), where the free end of the flag (in blue) is facing an incoming uniform water flow $U$.
The study of inverted flags has gained much interest in recent years \cite{flap1,flap2,flap3} for its potential applications in energy harvesting systems.

Across the parameter space of the inverted flag, including flag flexural rigidity and incoming flow properties, the inverted flag motion can display a variety of phenomena, including the undeformed equilibrium state, stable small-deflection state,  small-deflection flapping, large-amplitude periodic flapping and chaotic flapping  \cite{flap1}.
here we focus on the large-amplitude periodic flapping regime, wherein the inverted-flag system exhibits greater strain energy than conventional flag flapping \cite{flap1}.
These large bending strains make the inverted-flag system a promising candidate for energy harvesting technologies that convert strain energy to electricity.
In the large-amplitude flapping regime, the system has a stable limit cycle and three unstable fixed points, which have been identified numerically as steady-state solutions to fully coupled governing equations \cite{flap1}.
One fixed point is the undeformed equilibrium position, which becomes an unstable saddle point at sufficiently low bending stiffness.
The other two symmetrical fixed points are located between the saddle point and the limit cycle deflection amplitude.

In the experiments, an inverted flag with a white edge is released from different positions with near-zero velocity and allowed to evolve into a periodic orbit \cite{xu24}.
After preprocessing the videos with background removal, we track the endpoint of the inverted flag with our proposed algorithm.
Three frames selected at different times of a training video are shown in Figure \ref{fig:flag}(b).
To construct an SSM model, we train on four trajectories as shown in Figure \ref{fig:flag}(c).
The pixel coordinates have been normalized against the width of the video.

We delay-embed the $x$ coordinate signal of the flag end point to create a 5-dimensional observable space, in which \emph{SSMLearn} identifies a 2-dimensional, 3rd-order SSM.
Since the origin is a saddle point, \emph{SSMLearn} outputs reduced dynamics corresponding to two real eigenvalues with opposite sign.
We identify the reduced dynamics on the SSM up to 9th order, with the full equations presented in the Appendix.

In Figure \ref{fig:flag}(c), we apply the model to four unseen test trajectories released from different initial positions and plot the reconstructions with the full test dataset.
We find good qualitative agreement between the model predictions and the test data.
The discrepancies are likely caused by varying experimental conditions, such as room temperature and material fatigue, between videos \cite{xu24}.

In Figure \ref{fig:flag}(d), we construct a phase portrait by placing various initial conditions on the SSM and advect both forward and backward in time.
We find that we recover the stable (blue) and unstable (red) manifold of the saddle point and the three fixed points (green) of the system, along with the stable limit cycle.

\subsection{Liquid sloshing}
For our third example, we consider fluid oscillation in a partially filled container subject to horizontal external forcing.
The sloshing motion of the fluid can be strongly nonlinear \cite{sloshing}.
Increased fluid oscillation results in greater shearing at the tank wall, and the damping of the system grows nonlinearly with sloshing amplitude.
Further, the study in \cite{SSMLearn} found a softening response in the system, where the frequency decreases at a higher amplitude.
Understanding and modeling the system's nonlinear response is crucial for the design and analysis of structures that involve liquid-containing systems, such as fluid-transporting trucks \cite{truck_slosh} and spacecraft fuel tanks \cite{space_slosh}.

\begin{figure}[t]
    \centering
    \includegraphics[width=0.95\linewidth]{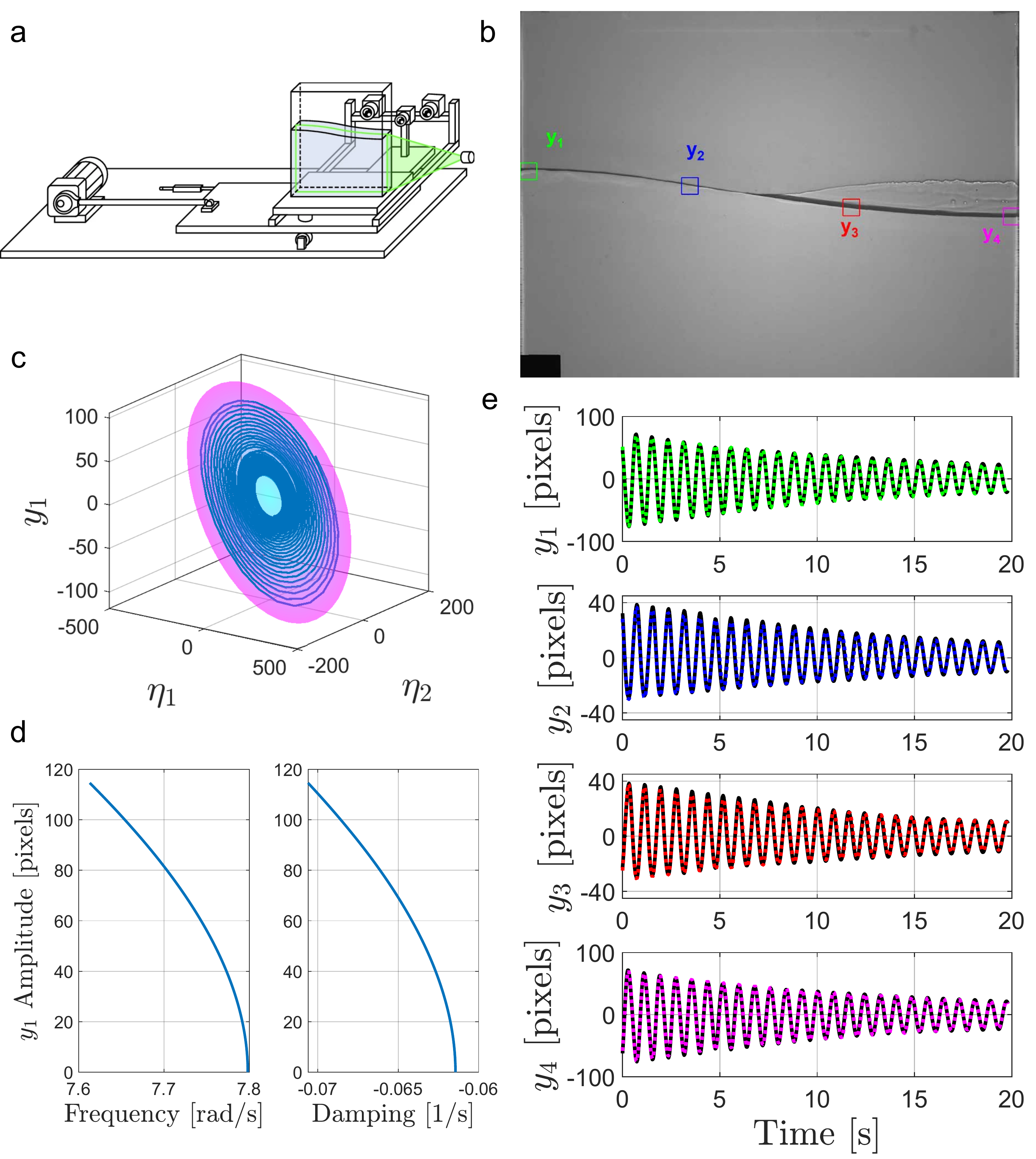}
    \caption[Slosh]%
    {\textbf{Data-driven Nonlinear Reduced-order Model on the Slowest SSM of Water Sloshing in a Tank} --- 
    (a) Sloshing experiment setup with partially filled tank and cameras, mounted on a platform linked to a motor \cite{sloshing}.
    (b) The four points $y_1,\dots,y_4$ at the water surface selected for tracking.
    (c) Training trajectory and identified SSM plotted in reduced coordinates and the leftmost vertical surface elevation $y_1$ in pixels.
    (d) Backbone and damping curves produced by the model.
    (e) Test trajectory for vertical positions of the four material points in black and the model predictions in colors corresponding to (b).
    }
    \label{fig:slosh}
\end{figure}

Here, we study videos from sloshing experiments performed by \cite{slosh_experiment}.
The experiments were performed in a rectangular tank of width 500 mm partially filled with water up to a height of 400 mm, as shown in Figure \ref{fig:slosh}(a).
The tank was mounted on a platform excited harmonically by a
motor.
Videos of the surface profile were recorded with a monochrome camera (1600x1200 pixels at 30 Hz) mounted on a frame moving with the tank.
The flow was illuminated by an LED light behind the tank, which creates a distinct shadow of the meniscus.
The light source was placed approximately 2 m away from the tank to provide uniform illumination and to avoid heating the fluid. 

In \cite{slosh_experiment}, the surface profile was extracted with a combination of image gradient and thresholding image processing techniques, specialized for the experimental setup and video format.
SSM-based models constructed from unforced decaying surface profile data were trained in \cite{SSMLearn,fastSSM}.
Here, we will apply our tracking method to extract the liquid surface heights at four, evenly spaced horizontal positions.
We initialize four templates, as shown in Figure \ref{fig:slosh}(b) with the search region constrained to only the vertical direction.
This is the only modification made to our tracking algorithm.
The horizontal constraint is required as the shadow of the meniscus is indistinguishable along the surface.
We assume the deformation of the surface is locally small, that is, the surface profiles within each template only undergo rigid rotations and vertical translations during sloshing motion.

In Figure \ref{fig:slosh}(e) we plot a decaying response signal extracted from experiment videos.
We use one such decaying trajectory for training and another one for testing.
We delay-embed the four extracted signals $y_{1,2,3,4}$ five times, with $\Delta t = 0.033 \mathrm{s}$, to create a 20-dimensional observable space, in which \emph{SSMLearn} identifies a 3rd-order, 2-dimensional SSM illustrated in Figure \ref{fig:slosh}(c).

On the SSM, we identify the reduced dynamics up to 3rd order and compute its 3rd order normal form
\begin{equation}
    \begin{aligned}
        \dot{\rho}\rho^{-1} &= -0.062 - 0.029\rho^2,\\
        \dot{\theta} &= 7.80 - 0.60\rho^2,
    \end{aligned}\label{eq:slosh_normal_form}
\end{equation}
where the first-order frequency term agrees with an analytical computation of the eigenfrequency \cite{slosh_experiment}.

In Figure \ref{fig:slosh}(d), we plot the nonlinear damping and softening response, with respect to the vertical pixels amplitude of the tracked point furthest to the left $y_1$, captured by our model.

In Figure \ref{fig:slosh}(e) we test our model on an unseen trajectory.
By inputting only the initial condition and integrating forward in time, our model reconstructs the full test trajectory within 4\% error.
We find our model to be accurate and correctly capture the essential nonlinearities of the system.

\subsection{Aerodynamic flutter}
As our fourth example, we consider aeroelastic flutter, a fluid-structure interaction phenomenon characterized by a self-excited structural oscillation wherein energy is drawn from the airstream through the movement of the structure \cite{intro_to_aeroelasticity}.
In a dynamical systems setting, flutter corresponds to a supercritical Hopf bifurcation where the bifurcation parameter is the airspeed.

\begin{figure}[ht]
    \centering
    \includegraphics[width=1\linewidth]{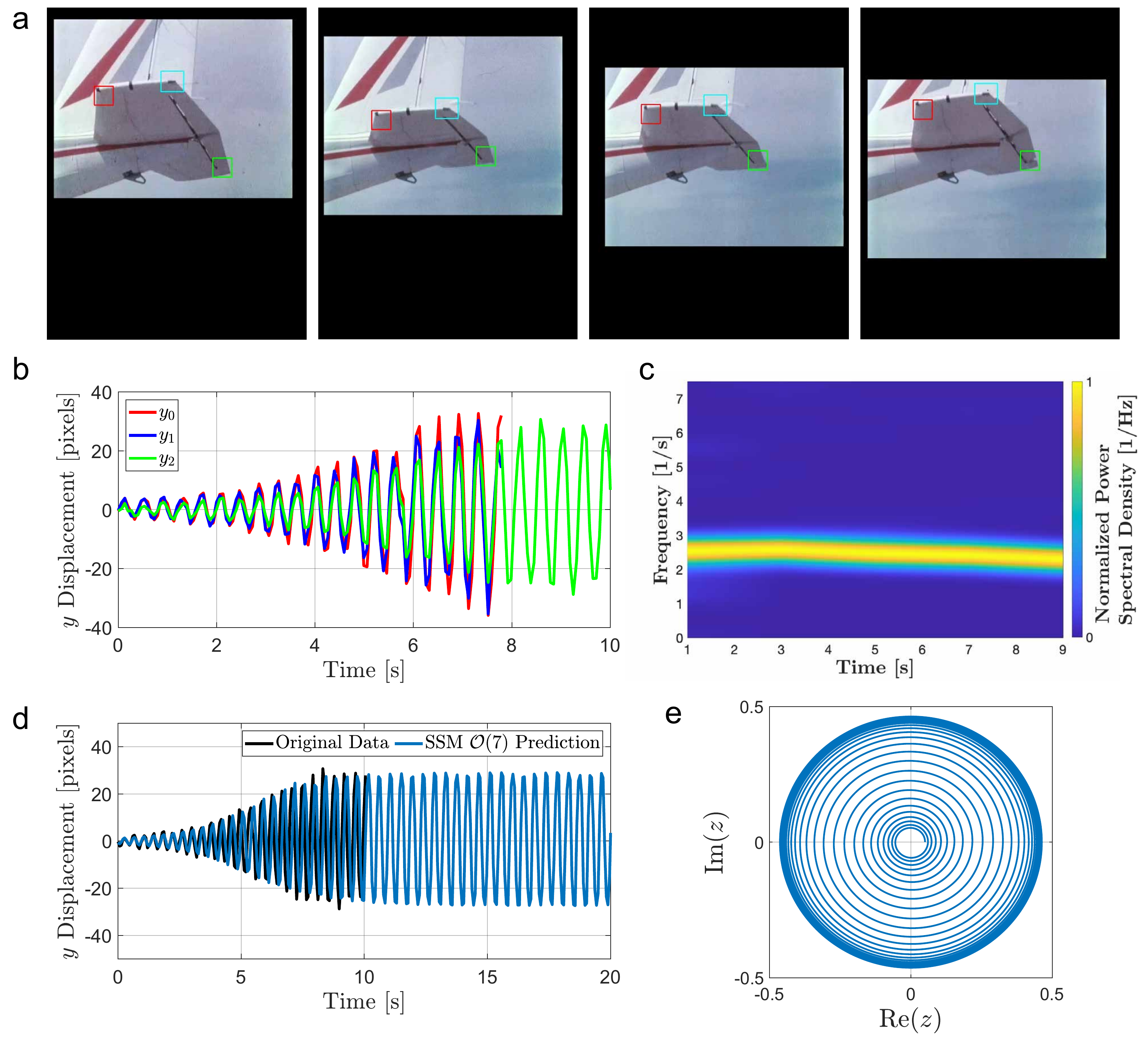}
    \caption[Flutter]%
    {\textbf{Data-driven Nonlinear Reduced-order Model on the Unstable Manifold for Tailplane Flutter} --- 
    (a) Four selected frames of the tailplane flutter video with the three tracked corners in different colors. The black margins of frames are paddings for frame stabilization since the subject is not static in the video.
    (b) Vertical displacements of the three tracked corners. Since the tailplane moves up and leaves the frame in the last two seconds of the video, there are no valid readings for the blue and red corners for $t>8$ s.
    (c) Spectrogram of the displacement signal at the bottom right corner (green) of the horizontal stabilizer.
    (d) Original trajectory of the vertical displacement of the bottom right corner, $y_2$, in black and its reconstruction in blue.
    By advecting the initial condition beyond the time of the original data, we observe a sustained and bounded limit cycle oscillation.
    (e) Reconstructed trajectory in the normal form coordinates.
    The model captures the limit cycle behavior.
    }
    \label{fig:flutter}
\end{figure}

Modeling aeroelastic flutter is important for the design of flexible aerodynamic structures.
In applications where the structural resilience of flexible bodies is crucial, the oscillatory motion plays a significant role in influencing the structure's dynamics and failure.
Aside from posing challenges in design scenarios, flutter also emerges as a method for harnessing energy \cite{flutter_energy}.
The simultaneous need to address and enhance flow-induced motion is thus important across diverse engineering fields.

Here, we consider a video of a flutter test from 1966 released by NASA \cite{flutter_video}.
The test was conducted with a Piper PA-30 Twin Commanche piloted by Fred Haise.
We assume the variation in the airspeed is small after the aeroelastic flutter initiation. 
As the horizontal tailplane in the video moves out of frame, we track the black hook at the bottom of the plane with our tracker and apply inverse translation to the video frame to apply frame stabilization.
Figure \ref{fig:flutter}(a) shows four stabilized frames of the video with three tracked corners of the tailplane in different colors.

We plot the extracted vertical displacements of the three corners in Figure \ref{fig:flutter}(b).
There is no reading for the blue and red corners in the last two seconds of the video as the plane moves partially out of the frame.
Hence we use only the bottom right corner of the tailplane (green bounding box and signal) to construct our reduced model.

For a two-degree-of-freedom (binary) aeroelastic model at moderate speed, coupling wing bending and torsion, the system has two complex conjugate pairs of eigenvalues corresponding to eigenmodes with positive damping and distinct frequencies \cite{intro_to_aeroelasticity}.
The flutter begins when the damping in one of the modes crosses zero, which occurs when the flow velocity exceeds a critical airspeed and the system undergoes a Hopf bifurcation.
The effective dampings in each of the motions must be simultaneously zero and the frequencies of both motions must be identical \cite{hancock1985teaching}.
Hence, the oscillation will have a single frequency.
This argument extends to more than two degrees of freedom.

Further, in Figure \ref{fig:flutter}(c), we plot a spectrogram of the training signal with the power spectral density normalized in each time slice. 
We find only a single dominant frequency is present in the signals, so we aim to identify a 2D SSM in an appropriate observable space for our model reduction.
We delay-embed the signal with $\Delta t = 0.0667 \mathrm{s}$, to construct an SSM-based model in a 5-dimensional observable space.
We identify and parameterize a 3rd-order, 2D SSM, on which \emph{SSMLearn} outputs a 7th-order reduced dynamics and computes a 7th-order normal form
\begin{align}
    \begin{split}
    \dot{\rho}\rho^{-1} &= +0.4844 - 1.679\rho^2 - 8.516\rho^4 + 27.28\rho^6,\\
    \dot{\theta} &= 15.90 - 34.64\rho^2 + 377.1\rho^4 - 1213\rho^6.
    \end{split}
\end{align}

In Figure \ref{fig:flutter}(d), we integrate in time the initial condition of the original trajectory beyond the original data, we observe the oscillation sustains for all future time and its amplitude remains bounded. 
The trajectory in the normal space is illustrated in Figure \ref{fig:flutter}(e), where the trajectory grows initially and enters a stable constant amplitude oscillation.
We find that the limit cycle is captured by our model and the dynamics predictions remain bounded.

\subsection{Wheel shimmy}
\begin{figure}[t]
    \centering
    \includegraphics[width=1\linewidth]{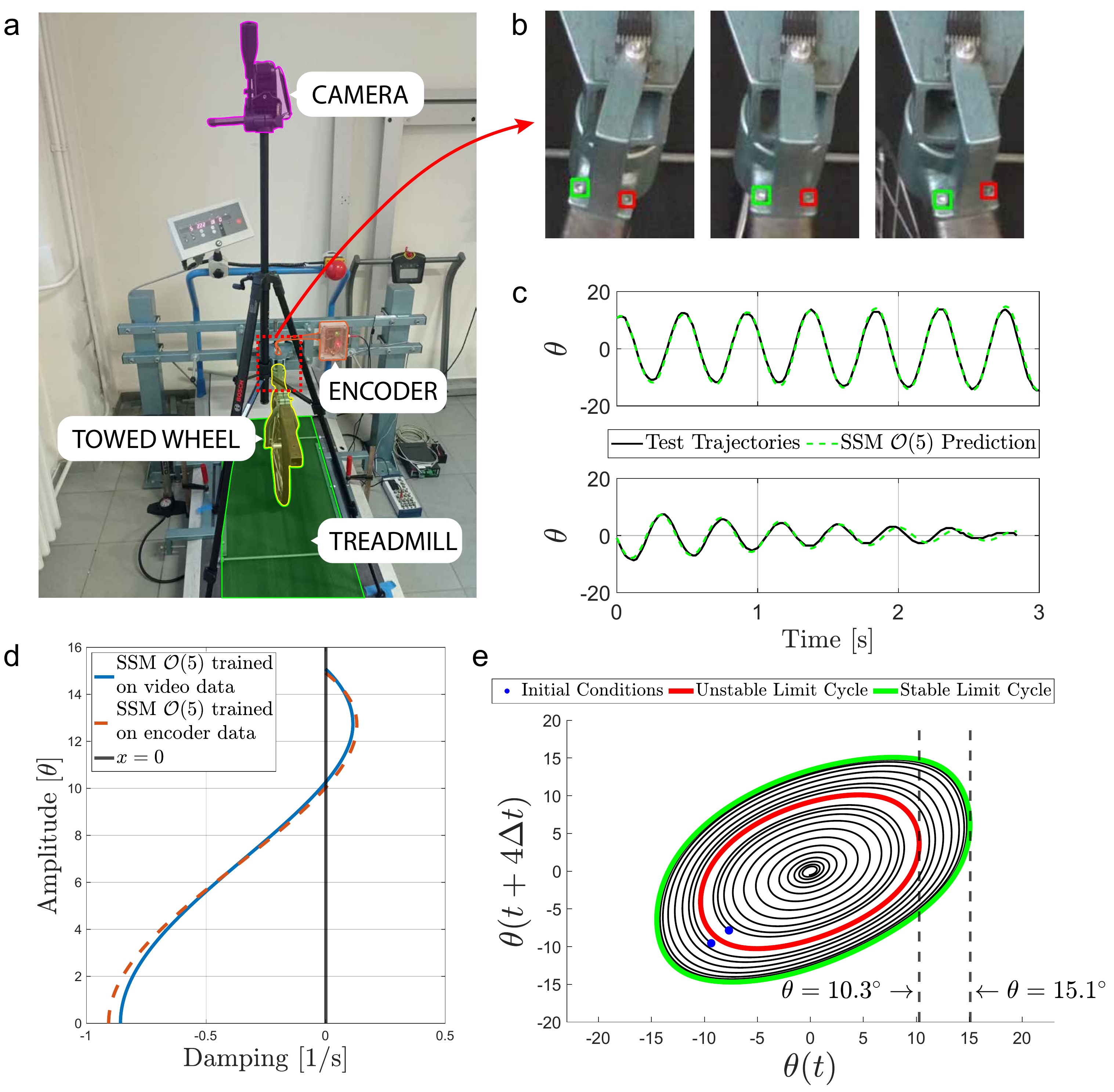}
    \caption[Shimmy]%
    {\textbf{Data-driven Nonlinear Reduced-order Model on the Slowest SSM of Wheel-Shimmy} --- 
    (a) The experimental rig with a treadmill in green, towed bicycle wheel in yellow, encoder in orange, and camera in magenta.
    The red box indicates the region we use for extracting data.
    (b) Three examples of processed frames, where the tracked bolts are highlighted with green and red bounding boxes.
    The $x$ and $y$ coordinates of the box centers are used to compute yaw angle $\theta$.
    (c) Test trajectories for limit cycle and decaying oscillations, and their reconstructions predicted by the trained model.
    (d) Damping curve output by the SSM-reduced model trained from video data and encoder data.
    The y-axis zero-crossings correspond to the unstable and stable limit cycle amplitudes.
    (e) Phase portrait of the system.
    }
    \label{fig:shimmy}
\end{figure}

As the last example, we consider the self-excited vibration of towed wheels, often referred to as wheel-shimmy.
This phenomenon can occur across vehicle systems, such as motorcycles \cite{shimmy_motorcycle}, tractors \cite{shimmy_tractor} and aircraft landing gears \cite{shimmy_landing_gear}, with significant safety implications.

The rich dynamics of the wheel-shimmy originate from the elastic tyres whose ground contact regions are subject to partial sticking and sliding \cite{beregi2019theoretical}.
The time delay in the tyre-ground contact and dry friction result in subcritical Hopf bifurcations in the infinite-dimensional and non-smooth system, giving rise to bistable parameter domains.
In the bistable parameter region, the stable rectilinear motion and periodic oscillation coexist with domains of attraction separated by unstable limit cycles.
The experimental, analytical, and numerical findings collectively affirm the presence of subcritical bifurcations in rectilinear motion, where bistable parameter domains were observed under a fixed caster length and a range of towing speeds \cite{beregi2019bifurcation, beregi2019analysis, beregi2019theoretical}.

Here, we study a caster-wheel system running on a treadmill with a fixed set of system parameters, caster length and tow speed, in a bistable domain.
We film top-down videos (1920x1080 px, 50 Hz) to generate training and test data.
The experimental setup is shown in Figure \ref{fig:shimmy}(a), which is a modification from the one studied in \cite{beregi2019bifurcation, beregi2019theoretical}.
On top of the treadmill, a bicycle wheel is fixed to a caster, which is mounted to the rack by a rotational joint.
A rotational encoder built on the joint provides yaw angle reading for validation.

We track two bolts (see Figure \ref{fig:shimmy}(b)) on the rotating caster bar to obtain $(x,y)$ pixel positions of the bolt centers, which we use to compute yaw angle readings.
Perturbations, by hand push, were applied to the wheel to generate trajectories that converge to the stable fixed point, $\theta = 0$, and the periodic orbit.
Two such trajectories are plotted in Figure \ref{fig:shimmy}(c).
We use a pair of trajectories, capturing the transients of decaying vibrations and large amplitude oscillations, for training and another pair for testing.

Frequency analysis shows that only a single frequency is present in the signals, so we aim to identify a 2D SSM in an observable space for our model reduction.
Based on Takens's embedding theorem, we delay-embed the signal five times at a time step $\Delta t = 0.02 \mathrm{s} = 1 / 50 \mathrm{Hz}$.
In the 5-dimensional observable space, \emph{SSMLearn} identifies a 3rd-order, 2-dimensional manifold, on which we identify the reduced dynamics in the 5th-order normal form
\begin{align}
    \begin{split}
    \dot{\rho}\rho^{-1} &= -0.8583 + 12.11\rho^2 - 37.71\rho^4,\\
    \dot{\theta} &= 15.17 - 9.155\rho^2 + 7.398\rho^4.
    \end{split}
\end{align}
Note that a 5th-order normal form is the minimal order to model the bistable dynamics of the system because we expect two zeros in the damping curve corresponding to the unstable and stable limit cycle amplitudes.

In Figure \ref{fig:shimmy}(c), we test our model by advecting two unseen initial conditions and comparing them against the full test trajectories.
With our training settings, we find a reconstruction error of 3.17\% for the limit cycle and 5.07\% for the decaying oscillation.

As a validation, we train a model using yaw angle readings from the encoder installed on the rotating joint (sampling rate of 200 Hz) with the same trajectories and training setting.
We plot the damping curve output by the two models, trained on video data and encoder data, in Figure \ref{fig:shimmy}(d).
The damping curves from the two models are consistent and predict a similar unstable limit cycle amplitude near $10^\circ$.
The small discrepancy near zero amplitude is caused by a higher signal-to-noise ratio from both measurement methods.
Figure \ref{fig:shimmy}(e) shows the phase portrait constructed by integrating two initial conditions near the unstable limit cycle.

\section{Conclusion}
We have developed a tracking algorithm based on rotation-aware template matching to extract data from experiment videos, which we have used to construct nonlinear reduced-order models with the open-source MATLAB package \emph{SSMLearn}. We have obtained SSM-based models that accurately captured nonlinear phenomena such as the nonlinear damping of a double pendulum, the softening response of water tank sloshing, multiple coexisting isolated steady states of an inverted flag, stable limit cycle oscillations in the tailplain flutter of an aircfraft and bistable dynamics in wheel shimmy.
This broad range of applications demonstrates the versatility of the tracking method and the general applicability of SSM-reduced models constructed from video data.

The main objective of conventional trackers is \emph{object localization}, whereas the purpose of our tracking algorithm is \emph{experimental sensing}.
Conventional filter-based trackers are unsuitable for this purpose because the tracked points or features are implicitly defined and updated at each frame, which leads to relatively higher noise, uncertainty in extracted time series data, and drifting issues (see Appendix for a comparison).

Although the tracking method devised here can extract accurate trajectories for the construction of reduced-order models, a limitation is its computational cost, which is notably higher than state-of-the-art CSRT trackers \cite{CSRT}, as well as its limited robustness to environmental changes.
Our tracking algorithm takes a brute-force approach to identify optimal matches in target frames, involving the computation of similarity scores between the template and all locations in the search region for multiple potential template rotations.
Even though the tracker utilizes information from previous matches to reduce the range of search regions and template rotations, the computational cost remains high.

In future work, we envision reducing the computational time significantly by parallelizing the matching computation without compromising its accuracy.
Additionally, one can improve computational efficiency by incorporating motion models to reduce search regions or by calculating matching correlations in the frequency domain similar to the technique used in DCF-based methods \cite{dcf_scale,CSRT}.
Exploring alternative handcrafted features like histograms of oriented gradients (HoG) or Gabor filters could also enhance tracking robustness but at the cost of reduced model interpretability.

Up to a rotation, our template matching tracking method assumes that variations in appearance and lighting are small and a single initialization of the template is therefore sufficient throughout the full video length.
It would be possible to extend our template matching to other classes of linear transformations to account for small general rotations and deformations. This, however,  would require high-fidelity video and improvements to the computational performance.

These limitations can be addressed in future developments of the algorithm presented here. Our general conclusion is that in  situations where the speed of extracting an SSM-reduced model is not critically important,  the present form our the algorithm already provides reliable models for nonlinear system identification, trajectory prediction, and localization of hidden features of the dynamics such as unstable limit cycles.

\backmatter
\section*{Declarations}

\subsection*{Funding}
This research was supported by the Swiss National Foundation (SNF) grant no. 200021-214908.

\subsection*{Conflict of interest}
The authors have no relevant financial or non-financial interests to disclose.

\subsection*{Author contributions}
A.Y. carried out the analysis and wrote the first draft of the paper; J.A. supervised the research and edited the paper; F.K. carried out the experimental measurements to the wheel shimmy example and edited the paper; G.S. supervised the experimental work and contributed to the analysis; G.H. designed the research, lead the research team and edited the paper.




\begin{appendices}
\section{Error metric}\label{app:ERMSE_CNMTE}

The root mean square error (RMSE) is a commonly used metric in statistics and machine learning to evaluate the accuracy or performance of a prediction or model. It is defined as
\begin{align*}
    \mathrm{RMSE} = \sqrt{\frac{1}{N} \sum_{i=1}^N (\hat{y}_i - y_i)^2},
\end{align*}
which would not suffice in our context where the observation vector $\bm{y}(t) \in \mathbb{R}^n$ is not a scalar and for each trajectory, each dimension of the $\bm{y}(t)$ could take different ranges of values.

In other words, the prediction error for each dimension of the output should be normalized by the corresponding range of the state. For instance, the error of the two states of a single pendulum, $\theta$ and $\dot{\theta}$, should be normalized separately as they could take different ranges of values: $\theta \in [0, 2\pi)$ and $\dot{\theta} \in (-\infty, \infty)$.

To quantify the performance of models that have multi-dimensional output, the author proposes to use a variation of RMSE extended to vector observations, which will be referred to as the extended root mean square error (ERMSE) defined as
\begin{align}\label{eq:ERMSE}
    \mathrm{ERMSE} =  \sqrt{\frac{1}{N}\sum_{i=1}^N \frac{1}{n} \sum_{j=1}^n \left( \frac{ \hat{y}_j(t_i) - y_j(t_i)}{ \max_i y_j(t_i) - \min_i y_j(t_i)} \right)^2},
\end{align}
for each $\bm{y}(t) \in \mathbb{R}^n$ output vector with $N$ samples and the corresponding reconstruction $\bm{\hat{y}}(t)$ from model.

Alternatively, a modified version of the root mean trajectory error (NMTE) defined in \cite{fastSSM} as
\begin{align*}
    \mathrm{NMTE} =  \frac{1}{N} \frac{1}{||\max_i{\bm{y}(t_i)}||}\sum_{i=1}^N ||\hat{\bm{y}}(t_i) - \bm{y}(t_i)||,
\end{align*}
can also be used. Instead of normalizing the average vectorial error against the maximum norm of $\bm{y}(t)$, the author suggests normalizing each dimension/component to its range separately and then computing the average. The author will refer to this metric as component-normalized mean trajectory error (CNMTE) defined as
\begin{align}\label{eq:CNMTE}
    \mathrm{CNMTE} =  \frac{1}{N}\sum_{i=1}^N\sqrt{ \sum_{j=1}^k \left( \frac{ \hat{y}_j(t_i) - y_j(t_i)}{ \max_i y_j(t_i) - \min_i y_j(t_i)} \right)^2}.
\end{align}

Even though the expression of (\ref{eq:ERMSE}) and (\ref{eq:CNMTE}) are similar, there is a difference in meaning. The ERMSE is an extension of RMSE, which aims to measure the performance of machine-learning models, where models take multiple inputs and produce multiple outputs of different scales. The outputs are not assumed to be related to one another. Whereas the CNMTE carries a physical meaning, where it takes the average of the norm of the normalized vectorial error of the states, which represents trajectory differences.

The choice of which metric to use thus depends on the setting. The author will use the ERMSE to quantify the error of manifold parameterization, as this step of the procedure involves fitting a polynomial to a set of scattered data points from all the training trajectories, and the CNMTE will be used to quantify the error of dynamics predictions from the reduced model, where the idea of comparing trajectories is important.

\section{Equations of motion for the double pendulum}\label{app:dp_derivation}

The derivation of the equations of motion begins with the definitions for the coordinates and velocity of the center of masses of the two pendulum rods
\begin{align*}
    (x_1, y_1) &= \left( \frac{1}{2} l_1 \sin{\theta_1}, \: - \frac{1}{2} l_1 \cos{\theta_1} \right), \\
    (\dot{x}_1, \dot{y}_1) &= \left( \frac{1}{2} l_1 \dot{\theta}_1 \cos{\theta_1}, \: \frac{1}{2} l_1 \dot{\theta}_1 \sin{\theta_1} \right), \\
    (x_2, y_2) &= \left( l_1 \sin{\theta_1} + \frac{1}{2} l_2 \sin{\theta_2}, \: - l_1 \cos{\theta_1} - \frac{1}{2} l_2 \cos{\theta_2} \right), \\
    (\dot{x}_2, \dot{y}_2) &= \left( l_1 \dot{\theta}_1 \cos{\theta_1} + \frac{1}{2} l_2 \dot{\theta}_2 \cos{\theta_2}, \: l_1 \dot{\theta}_1 \sin{\theta_1} + \frac{1}{2} l_2 \dot{\theta}_2 \sin{\theta_2} \right),
\end{align*}
where $g$ is the gravitational acceleration constant. 
These definitions assume the center of masses of the two rods coincides with the halfway point of the physical rod lengths.

The kinetic energy $T$ and potential energy $V$ are
\begin{align*}
    T &= \frac{1}{2} m_1(\dot{x}_1^2 + \dot{y}_1^2) + \frac{1}{2} m_2(\dot{x}_2^2 + \dot{y}_2^2) + \frac{1}{2} I_1 \dot{\theta}_1^2 + \frac{1}{2} I_2 \dot{\theta}_2^2,\\
    \begin{split}
      &= \underbrace{\left(\frac{1}{2}m_1 \left(\frac{1}{2}l_1\right)^2 + \frac{1}{2}I_1^2 + \frac{1}{2}m_2 l_1^2\right)}_{A} \dot{\theta}_1^2 + \underbrace{\left(\frac{1}{2}m_2 \left(\frac{1}{2}l_2\right)^2 + \frac{1}{2}I_2^2\right)}_{B} \dot{\theta}_2^2 \\
      &+ \underbrace{\frac{1}{2} m_2 l_1 l_2}_{C} \dot{\theta}_1 \dot{\theta}_2 \cos{(\theta_2 - \theta_1)},
    \end{split}\\
    V &= m_1 g y_1 + m_2 g y_2, \\
      &= - \underbrace{g l_1 \left(\frac{1}{2} m_1 + m_2 \right)}_{D} \cos{\theta_1} - \underbrace{\frac{1}{2} m_2 g l_2}_{E} \cos{\theta_2},
\end{align*}
where $I_{1,2}$ are the associated moment of inertia, and $A, B, C, D, E$ are all constants consisting of system parameters.

Invoking the parallel axis theorem, we compute the moment of inertia of the pendulum as a summation of contribution from a rectangular prism and two semicylinders at both ends.
\begin{align*}
    I_i = \frac{1}{12}M_{rod}(l_i^2 + w_i^2) + 2M_{semicylinder}\left[\left(\frac{1}{16}-\frac{4}{9\pi^2}\right)(w_i)^2 + \left(\frac{l_i}{2}+\frac{2w_i}{3\pi}\right)^2\right],
\end{align*}
where $w_i$ corresponds to the width of the pendulum arms.

Thus, we can immediately write the Lagrangian of the conservative system as
\begin{align*}
\mathcal{L} = T - V.
\end{align*}

The equations of motion of the conservative system could be obtained by solving the Euler-Lagrange Equations
\begin{align}\label{eq:conserve_euler_lagrange}
    \frac{\mathrm{d}}{\mathrm{d}t} \left(\frac{\partial \mathcal{L}}{\partial \dot{q_i}} \right) - \frac{\partial \mathcal{L}}{\partial q_i} = 0, \; \mathrm{for} \, i = 1, 2,
\end{align}
which yield
\begin{align}
    \begin{split}\label{eq:dp_conserve1}
        \ddot{\theta}_1 = \frac{1}{K} \Biggr[C^2\sin{(\theta_1 -\theta_2)}\cos{(\theta_1 -\theta_2)}\dot{\theta}_1^2  + 2BC\sin{(\theta_1 - \theta_2)}\dot{\theta}_2^2 \\
        + 2BD\sin{\theta_1} - CE\cos{(\theta_1-\theta_2)}\sin{\theta_2} \Biggr],
    \end{split}
\end{align}
\begin{align}
        \begin{split}\label{eq:dp_conserve2}
        \ddot{\theta}_2 = \frac{1}{K} \Biggr[-2AC\sin{(\theta_1 -\theta_2)}\dot{\theta}_1^2 - C^2\sin{(\theta_1-\theta_2)}\cos{(\theta_1-\theta_2)}\dot{\theta}_2^2\\ 
        - CD\cos{(\theta_1-\theta_2)}\sin{\theta_1} + 2AE\sin{\theta_2} \Biggr],
    \end{split}
\end{align}
\begin{align*}
    \mathrm{where}\;
    \begin{cases}
        &A = \frac{1}{2}m_1 \left(\frac{1}{2}l_1\right)^2 + \frac{1}{2}I_1^2 + \frac{1}{2}m_2 l_1^2,\\
        &B = \frac{1}{2}m_2 \left(\frac{1}{2}l_2\right)^2 + \frac{1}{2}I_2^2,\\
        &C = \frac{1}{2}m_2 l_1 l_2,\\
        &D = (\frac{1}{2}m_1 + m_2)gl_1,\\
        &E = \frac{1}{2}m_2 g l_2,\\
        &K = C^2 \cos^2{(\theta_1-\theta_2)} - 4AB.
    \end{cases} 
\end{align*}

However, without dissipation, such a system is chaotic and would not suffice as a good model of the real system in this study. Hence, a dissipation must be included in the analytical model to replicate the physical system.

\subsection{Rayleigh dissipation function}
In 1881, Lord Rayleigh demonstrated that when a dissipative force $\bm{F}$ is proportional to velocity, it can be represented by a scalar potential that depends on the generalized velocities $\bm{\dot{q}}$. This scalar potential is known as the Rayleigh dissipation function $\mathcal{D}$
\begin{align}\label{eq:rayleigh_function}
    \mathcal{D}(\bm{\dot{q}}) \equiv \frac{1}{2}\sum_{i=1}^n\sum_{j=1}^n \beta_{ij}\dot{q_i}\dot{q_j},
\end{align}
and it provides an elegant way to incorporate linear velocity-dependent dissipative forces in Lagrangian and Hamiltonian mechanics.

Dissipative drag force in the direction of $\bm{q}$ is then defined as the negative velocity gradient of the dissipation function (\ref{eq:rayleigh_function}), that is
\begin{align*}
    \bm{F} = -\nabla_{\bm{\dot{q}}}\mathcal{D}(\bm{\dot{q}}).
\end{align*}

In the context of a double pendulum, it is reasonable to assume that frictions from the two joints are the main sources of dissipation in the system as opposed to other factors such as air resistance.

If we consider the two joints as free bodies, the dissipation from the first joint only depends on the rate of rotation of the upper arm, whereas the second joint linking the two arms depends on the relative rate of rotation of the two rods. Thus, the dissipation function for the system is of the form
\begin{align*}
    \mathcal{D} =  \underbrace{\frac{\beta_1}{2}\dot{\theta}_1^2}_{\mathrm{Joint 1}} + \underbrace{\frac{\beta_2}{2}(\dot{\theta}_2 - \dot{\theta}_1)^2}_{\mathrm{Joint 2}} \implies
    \begin{cases}
        F_1 &= - \frac{\partial \mathcal{D}}{\partial \dot{\theta}_1} = -(\beta_1 + \beta_2)\dot{\theta}_1 + \beta_2\dot{\theta}_2\\
        F_2 &= - \frac{\partial \mathcal{D}}{\partial \dot{\theta}_2} = \beta_2\dot{\theta}_1 -\beta_2\dot{\theta}_2
    \end{cases}.\nonumber
\end{align*}

The dissipative generalized forces can then be added to the Euler-Lagrange equation (\ref{eq:conserve_euler_lagrange}) to obtain
\begin{align}\label{eq:dissipate_euler_lagrange}
    \frac{\mathrm{d}}{\mathrm{d}t} \left(\frac{\partial \mathcal{L}}{\partial \dot{q_i}} \right) - \frac{\partial \mathcal{L}}{\partial q_i} = F_i, \; \mathrm{for} \; i = 1,2.
\end{align}

Solving (\ref{eq:dissipate_euler_lagrange}) arrives at a modified version of (\ref{eq:dp_conserve1}) and (\ref{eq:dp_conserve2})
\begin{align}
    &\begin{split}
        \ddot{\theta}_1 &= \frac{1}{K} \left[... + 2B\left((\beta_1 + \beta_2)\dot{\theta}_1 - \beta_2\dot{\theta}_2\right) + C\cos{(\theta_1 - \theta_2)} \left(\beta_2\dot{\theta}_1 - \beta_2\dot{\theta}_2 \right) \right]\\
        &= f_3(\theta_1, \theta_2, \dot{\theta}_1, \dot{\theta}_2, \bm{\mu}), \label{eq:dp_full1}\\
    \end{split}\\
    &\begin{split}
        \ddot{\theta}_2 &= \frac{1}{K} \left[... - 2A\left(\beta_2\dot{\theta}_1 - \beta_2\dot{\theta}_2 \right) - C\cos{(\theta_1 - \theta_2)}\left((\beta_1 + \beta_2)\dot{\theta}_1 - \beta_2\dot{\theta}_2\right) \right]\\
        &= f_4(\theta_1, \theta_2, \dot{\theta}_1, \dot{\theta}_2, \bm{\mu}), \label{eq:dp_full2}
    \end{split}
\end{align}
where $\bm{\mu}$ is a vector collecting all system parameters. This is the equation of motion for a double pendulum with modelled joint dissipation.
Even though the damping model linearly depends on angular velocities, it became nonlinear due to geometric nonlinearities.

\section{SSM-reduced model for the inverted flag}
\begin{align*}
    &\xi_1 = \\
    &- 2.325\xi_2 + 0.3512\xi_2^2 + 1.132\xi_2^3 - 1.479\xi_2^4 - 1.806\xi_2^5 + 1.806\xi_2^6\\
    &+ 1.742\xi_2^7 - 0.6762\xi_2^8 - 0.5968\xi_2^9 + 0.09106\xi_1 + 0.128\xi_1\xi_2 - 0.8056\xi_1\xi_2^2\\
    &- 0.5981\xi_1\xi_2^3 + 5.407\xi_1\xi_2^4 + 0.3837\xi_1\xi_2^5 - 8.837\xi_1\xi_2^6 + 0.1211\xi_1\xi_2^7 + 4.119\xi_1\xi_2^8\\
    &+ 0.2219\xi_1^2 - 0.502\xi_1^2\xi_2 - 1.619\xi_1^2\xi_2^2 + 0.5419\xi_1^2\xi_2^3 + 3.104\xi_1^2\xi_2^4 + 0.7389\xi_1^2\xi_2^5\\
    &- 1.701\xi_1^2\xi_2^6 - 0.608\xi_1^2\xi_2^7 + 0.6853\xi_1^3 + 0.3954\xi_1^3\xi_2 - 3.645\xi_1^3\xi_2^2 + 0.1459\xi_1^3\xi_2^3\\
    &+ 2.215\xi_1^3\xi_2^4 - 0.3787\xi_1^3\xi_2^5 + 0.6783\xi_1^3\xi_2^6 - 0.539\xi_1^4 + 3.49\xi_1^4\xi_2 + 1.664\xi_1^4\xi_2^2\\
    &- 5.234\xi_1^4\xi_2^3 - 1.183\xi_1^4\xi_2^4 + 1.699\xi_1^4\xi_2^5 - 1.499\xi_1^5 - 0.8057\xi_1^5\xi_2 + 5.23\xi_1^5\xi_2^2\\
    &+ 0.4301\xi_1^5\xi_2^3 - 2.93\xi_1^5\xi_2^4 + 0.4203\xi_1^6 - 3.089\xi_1^6\xi_2 - 0.5586\xi_1^6\xi_2^2 + 2.61\xi_1^6\xi_2^3\\
    &+ 0.9821\xi_1^7 + 0.2868\xi_1^7\xi_2 - 1.658\xi_1^7\xi_2^2 - 0.1005\xi_1^8 + 0.7641\xi_1^8\xi_2 - 0.2007\xi_1^9
    \\
    &\xi_2 = \\
    &- 0.2462\xi_2 + 0.6229\xi_2^2 + 5.107\xi_2^3 - 1.724\xi_2^4 - 16.8\xi_2^5 + 1.803\xi_2^6\\
    &+ 19.39\xi_2^7 - 0.6994\xi_2^8 - 7.426\xi_2^9 - 2.741\xi_1 + 1.084\xi_1\xi_2 + 11.86\xi_1\xi_2^2\\
    &- 5.412\xi_1\xi_2^3 - 10.74\xi_1\xi_2^4 + 6.211\xi_1\xi_2^5 - 1.404\xi_1\xi_2^6 - 1.81\xi_1\xi_2^7 + 3.985\xi_1\xi_2^8\\
    &+ 0.1505\xi_1^2 + 0.3348\xi_1^2\xi_2 - 4.276\xi_1^2\xi_2^2 - 8.45\xi_1^2\xi_2^3 + 8.499\xi_1^2\xi_2^4 + 16.54\xi_1^2\xi_2^5\\
    &- 4.433\xi_1^2\xi_2^6 - 8.455\xi_1^2\xi_2^7 + 8.555\xi_1^3 - 0.09146\xi_1^3\xi_2 - 43.53\xi_1^3\xi_2^2 + 4.224\xi_1^3\xi_2^3\\
    &+ 58.71\xi_1^3\xi_2^4 - 3.911\xi_1^3\xi_2^5 - 23.09\xi_1^3\xi_2^6 - 0.639\xi_1^4 + 7.441\xi_1^4\xi_2 + 4.639\xi_1^4\xi_2^2\\
    &- 10.29\xi_1^4\xi_2^3 - 4.067\xi_1^4\xi_2^4 + 2.415\xi_1^4\xi_2^5 - 10.49\xi_1^5 - 1.479\xi_1^5\xi_2 + 39.89\xi_1^5\xi_2^2\\
    &+ 0.04915\xi_1^5\xi_2^3 - 28.55\xi_1^5\xi_2^4 + 0.6601\xi_1^6 - 8.282\xi_1^6\xi_2 - 1.522\xi_1^6\xi_2^2 + 7.384\xi_1^6\xi_2^3 \\
    &+ 5.691\xi_1^7 + 0.61\xi_1^7\xi_2 - 10.51\xi_1^7\xi_2^2 - 0.1828\xi_1^8 + 2.22\xi_1^8\xi_2 - 1.073\xi_1^9
\end{align*}

\section{Comparison with other correlation filter-based trackers}
\begin{figure}[h]
    \centering
    \includegraphics[width=1\linewidth]{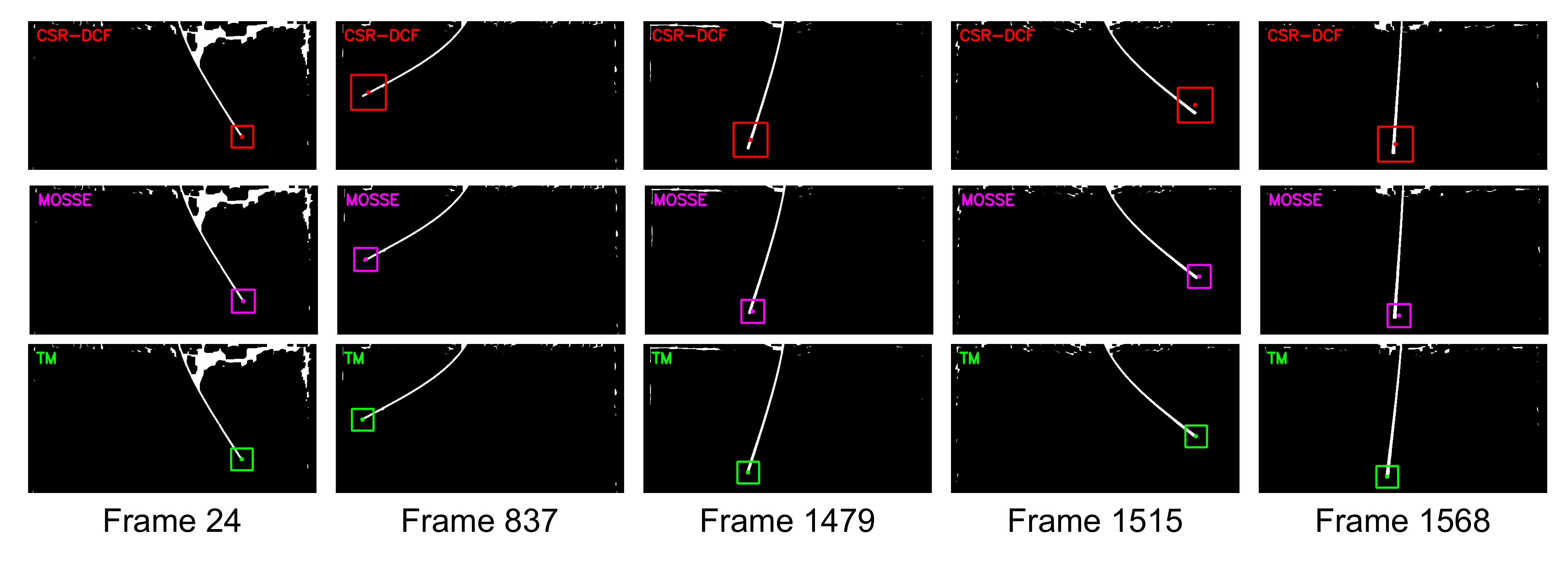}
    \caption[Comparison of Proposed Template Matching (TM) Tracker with CSR-DCF and MOSSE Tracker]%
    {\textbf{Comparison of Proposed Template Matching (TM) Tracker with CSR-DCF and MOSSE Tracker} --- This figure illustrates the qualitative performance comparison in an inverted flag experiment video.
    The proposed TM tracker, along with two advanced correlation filter-based trackers (CSR-DCF \cite{CSRT} and MOSSE \cite{MOSSE}), were initialized using an identical template in frame 24, followed by tracking the same video sequence.
    The center of the bounding box corresponds to the tracked point.
    Five frames showcase the drift issues common in traditional trackers.
    The first row (in red) are results from the CSR-DCF tracker.
    The second row (in magenta) are from the MOSSE tracker.
    The bottom row (in green) are from our proposed TM tracker.}
    \label{fig:compare_trackers}
\end{figure}

In Figure \ref{fig:compare_trackers}, we compare the tracking results of the proposed tracking method and two state-of-the-art trackers: CSR-DCF and MOSSE.
We initialize all three trackers with the same template with bounding box centers corresponding to the point we want to track, in this case, the tip of the inverted flag.
This point acts as the reference point for a qualitative evaluation of tracking accuracy.
We find our tracking algorithm to be accurate and does not exhibit any drifting issues as observed by the application of the other two methods.



\end{appendices}


\bibliography{sn-bibliography}

\end{document}